\documentclass[12pt,twoside]{amsart}

\usepackage{fullpage}

\usepackage{amsfonts}
\usepackage{amssymb}
\usepackage{amsmath}
\usepackage{epsfig}
\usepackage{psfrag}

\newtheorem{theorem}{Theorem}[section]
\newtheorem{prop}[theorem]{Proposition}
\newtheorem{lem}[theorem]{Lemma}
\newtheorem{cor}[theorem]{Corollary}

\numberwithin{equation}{section}

\def \lKN{\overset{\leftarrow}{K}^{\underset{N}{}}}

\def \lMN{\overset{\leftarrow}{M}^{\underset{N}{}}}

\def \R{\mathbb R}

\def \N{\mathbb N}

\def \E{\mathbb E}

\def \Tr{\mathcal{T}}
\def \Km{K^-}
\def \Kp{K^+}

\def \div{\mathrm {div}}

\long\def\symbolfootnote[#1]#2{\begingroup%
\def\thefootnote{\fnsymbol{footnote}}\footnote[#1]{#2}\endgroup}

\begin{document}

\title{Probabilistic analysis of the upwind scheme for transport}


\maketitle

\begin{center}
{\sc Fran\c{c}ois Delarue$^{1,2}$ and 
Fr\'ed\'eric Lagouti\`ere$^{1,3}$}
\end{center}

\symbolfootnote[0]{$^1$ Universit\'e Paris Diderot-Paris 7. 
E-mail: delarue@math.jussieu.fr,
lagoutie@math.jussieu.fr} 

\symbolfootnote[0]{$^2$ CNRS, UMR 7599, Laboratoire de Probabilit\'e et
Mod\`eles Al\'eatoires, F-75252, Paris, France. }

\symbolfootnote[0]{$^3$ CNRS, UMR 7598, Laboratoire Jacques-Louis
Lions, F-75005, Paris, France. }

\begin{abstract}
We provide a probabilistic analysis of the upwind scheme for
$d$-dimensional transport
equations. We associate a Markov chain with the numerical scheme and
then obtain a backward representation formula of Kolmogorov type for
the numerical solution. We then understand that the error induced by
the scheme is governed by the fluctuations of the Markov chain around
the characteristics of the flow. We show, in various
situations, that the fluctuations are of diffusive type. As a
by-product, we recover recent results due to Merlet and Vovelle
\cite{merlet:vovelle} and Merlet \cite{Mer?}: we prove that the scheme
is of order $1/2$ in $L^{\infty}([0,T],L^1(\R^d))$ for an initial
datum in $BV(\R^d)$ and of order $1/2-\varepsilon$, for all
$\varepsilon >0$, in $L^{\infty}([0,T] \times \R^d)$ for an initial
datum in $W^{1,\infty}(\R^d)$. Our analysis provides a new
interpretation of the {\em numerical diffusion} phenomenon. \\

\noindent
{\sc R\'esum\'e.}
Nous proposons une analyse probabiliste du sch\'ema upwind pour les
\'equations de transport en dimension $d$ quelconque. Pour cela, nous
associons au sch\'ema une cha\^ine de Markov qui nous permet d'obtenir
une formule de repr\'esentation de type Kolmogorov pour la solution
num\'erique. Nous comprenons alors que l'erreur due au sch\'ema est
gouvern\'ee  par les fluctuations de la cha\^ine de Markov autour des
caract\'eristiques du transport. Nous montrons, dans des situations
diverses, que ces fluctuations sont de type diffusif. Comme
cons\'equence, nous retrouvons des r\'esultats r\'ecents de Merlet et
Vovelle \cite{merlet:vovelle} et Merlet \cite{Mer?}~: nous montrons
que le sch\'ema upwind est d'ordre $1/2$ dans
$L^{\infty}([0,T],L^1(\R^d))$ pour une donn\'ee initiale dans
$BV(\R^d)$, et d'ordre $1/2 - \varepsilon$ pour tout $\varepsilon > 0$
pour une donn\'ee initiale dans $W^{1,\infty}(\R^d)$. Cette analyse
donne une interpr\'etation nouvelle du ph\'enom\`ene de {\em diffusion
num\'erique}. 
\\


\end{abstract}
\section{Introduction}

This paper provides a new analysis of the {\em upwind scheme} for the
transport problem in dimension $d \in \N \setminus \{ 0 \}$
\begin{equation}\label{trpb}
\left\{\begin{array}{l}
\partial_t u(t,x) + \left< a(x) , \nabla u (t,x) \right> = 0, \ (t,x)
\in [0, T] \times \R^d, \\
u(0,x) = u^0(x), \ x \in \R^d.
\end{array}\right.
\end{equation}
We assume $a$ to be Lipschitz continuous, so that \eqref{trpb} is
well-posed. Several different regularity assumptions are made on $u^0$
in the following, among which $u^0 \in W^{1,\infty}(\R^d)$ and $u^0
\in BV(\R^d)$. In any case, the unique solution to
\eqref{trpb} is $u(t,x) = u^0(Z(t,x))$ where $Z$ is the backward
characteristic, i.e. the solution of
\begin{equation}\label{carac}
\left\{\begin{array}{l}
\partial_t Z(t,x) = -a(x), \ (t,x) \in [0, T] \times \R^d, \\
Z(0,x) = x, \ x\in \R^d.
\end{array}\right.
\end{equation}

The upwind scheme is a standard method to solve this problem in an
approximate way (see for instance \cite{EymGalHer00}). It is derived
and described in Section \ref{upwind} below.

In dimension 1, the scheme is known to be first order consistent (with
respect to the maximal cell diameter $h$) with the transport
equation. Thus, it is first order convergent, for any $u^0 \in
\mathcal{C}^2(\R^d)$, provided that a Courant-Friedrichs-Lewy (CFL)
stability condition holds. For non-smooth initial data, the upwind
scheme is just of order 1/2. This loss of convergence order is
traditionally attributed to the {\em dissipative} character of the
scheme. Up to now, the use of the word ``dissipative'' has been
justified by the following fact: on a uniform mesh, the scheme is
second order consistent with an advection-diffusion equation, the
diffusion coefficient being first order with respect to $h$; as a
consequence, the numerical error at time $t$ is proportional to
$\sqrt{t h}$ for a non-smooth initial datum.

In this paper, we provide another explanation of the diffusive
behavior, which is valid on any general mesh in dimension $d$. We here
interpret the numerical diffusion by means of a stochastic
process. Let us briefly describe the basic idea. The approximate value
given by the upwind scheme in a cell $K$ at time step $n+1$ is a
convex combination of the approximate values at time step $n$ in the
cells neighboring $K$: 
$$
u_K^{n+1} = \sum_{L \in \Tr} p_{K,L} u_L^n,
$$
where $\Tr$ is the set of cells, and $p_{K,L} \in [0, 1]$ with
$\sum_{L \in \Tr} p_{K,L} = 1$ (see Section \ref{upwind} for the
complete definition of the scheme). This convex combination allows a
probabilistic interpretation: we can define a random sequence of cells
$\left( K_n \right)_{n\in\N}$ as a Markov chain with probability
transition, from $K$ to $L$, $p_{K,L}$. In this framework, the upwind
scheme appears as the {\em expectation of a random scheme} associated
with the chain $(K_n)_{n \geq 0}$. Precisely, the value $u_K^n$ is the
expectation of the value of $u^0$ in the cell $K_n$ when $K$ is chosen
as the starting cell of the chain. In a probabilistic way, we write:
$$
u_K^n = \E_K \left(u_{K_n}^0 \right),
$$
the symbol $K$ in the notation $\E_K$ meaning that $K_0=K$. (See
Theorems \ref{th:feynman_kac} and \ref{th:d2_feynman}.) In the theory
of stochastic processes, the above identity is a {\em backward
Kolmogorov formula}: it is the analogue of the representation formula
of the heat equation by the Brownian motion. We then understand the
chain $(K_n)_{n \geq 0}$ as a random backward characteristic.

Our main idea consists in analyzing the behavior of the random
characteristic according to the following program. The first point
is to show that the mean trend of the random characteristic coincides
with the exact characteristic $Z$, solution to \eqref{carac}. The
next step is to understand that the error of the numerical scheme is
governed by the fluctuations of the random characteristic around the
exact one. Heuristically, the order of the fluctuations is given by
the central limit theorem: therefore, we expect them to be controlled,
in a suitable sense, by $C h^{1/2}$ where $C$ only depends on the
datum $a$ and the time $t$. The final step is to derive the 1/2 order
of the scheme.

Applying this program, we establish the $1/2$ order in
$L^{\infty}([0,T],L^1({\mathbb R}^d))$ for $u^0 \in BV({\mathbb
R}^d)$. (See Theorem \ref{th:L_1BV}.) For $u^0 \in W^{1,\infty}({\mathbb
R}^d)$, we also prove that the scheme is of order $1/2-\varepsilon$
in $L^{\infty}([0,T],L^{\infty}({\mathbb R}^d))$ for all
$\varepsilon >0$. (See Theorem \ref{prop:L_inf}.) In this last case,
it is clear that $1/2$ is an upper bound (see \cite{God79}, and
\cite{TanTen95} for the non-linear case), but the exact convergence
order remains unknown. The reason why our estimate is better in the
$L^1$-in space norm 
may be explained as follows. Estimating the error in
$L^1$ amounts to average the initial cell of the random
characteristic. This additional averaging reduces the weight
of the trajectories of the chain that are away from $Z$.

Since the pioneering article of Kuznetsov \cite{Kuz76}, in which the
1/2 order is established in the Cartesian framework (for linear an
non-linear scalar equations), many papers have dealt with the rate
of convergence of the upwind scheme. Let us briefly review them.

For general scalar equations with a datum $u^0 \in
BV(\R^d)$, Cockburn, Coquel and Le Floch \cite{CocCoqLeF94}, Vila
\cite{Vil94}, \cite{BouPer98} and Chainais-Hillairet \cite{Cha99}
prove the (non-optimal) 1/4 order in $L^\infty([0, T], L^1(\R^d))$
under slightly different hypotheses (and for several schemes,
including the upwind one). 

For hyperbolic Friedrichs systems, Vila and Villedieu
\cite{VilVil03} derive a 1/2 order estimate in the $L^2([0, T]
\times \R^d)$ norm for $H^1(\R^d)$ initial data.

In the frame of the linear transport equation, which we are involved
in, Despr\'es \cite{Des04} proves a $1/2$ order estimate in the
$L^\infty([0, T], L^2(\R^d))$ norm in the case of $H^2(\R^d)$ data.
His proof relies on a precise study of the consistency of the scheme
after several time steps. (It is indeed known that the scheme is not
consistent at each time step on a general mesh.) 
For $\mathcal{C}^2$ initial data, Bouche, Ghidaglia and Pascal
\cite{BouGhiPas05} show the order 1 in the $L^\infty$ norm, under a
condition on the mesh that is related to consistency.
At last, in recent works,
\begin{itemize}
\item for an initial datum in $BV(\R^d)$, Merlet and Vovelle
\cite{merlet:vovelle} show the optimal estimate of order 1/2  in the
$L^\infty([0, T], L^1(\R^d))$ norm, 
\item for a Lipschitz continuous initial datum, Merlet \cite{Mer?}
shows the order $1/2 - \varepsilon$, for any $\varepsilon > 0$, in the
$L^\infty$ norm.
\end{itemize}

It is thus understood that our paper provides a new proof of the
results obtained in \cite{merlet:vovelle} and \cite{Mer?}. Actually,
our framework is slightly different since we do not assume the
velocity $a$ to be divergence-free, as done therein, but we assume it
to be independent of time. 
We think that this does not make fundamental
differences. Despite the similarity of our results, we insist on the
fact that the arguments here are completely different. As said
above, our proofs rely on the analysis of the stochastic
characteristic $(K_n)_{n \geq 0}$ (that shall mimic the exact
characteristic $Z$). In particular, we do not use energy estimates.
(Except those of Despr\'es and Bouche {\em et al.} based on the
consistency of the scheme, all the papers mentionned above are built
on energy or entropy estimates).

Our paper is organized as follows. In Section \ref{upwind}, we state
the framework of our analysis. In Section \ref{d1}, we focus on the
one-dimensional case to introduce, with great care, the notion of
stochastic characteristic. By the way, we establish a refined
estimate of the order of the scheme in the specific case where the
velocity is constant and the mesh is regular.
(See Proposition \ref{prop:TLC}). In Section \ref{sec:2}, we extend
the probabilistic interpretation of the upwind scheme to the higher
dimensional setting. We then provide a direct proof of the 1/2 order
in $L^\infty([0, T], L^1(\R^d))$ in the following simple case: $u^0$
is assumed to be periodic, as well as the mesh, and Lipschitz
continuous. This section is the heart of the paper. Refining the
strategy, we finally obtain in Section \ref{gensec} the announced
results. This last part is a bit more technical and relies on
a concentration inequality for martingales, which is given in
Annex, see Section \ref{annex}. 

\section{Framework and useful notations}\label{upwind}

Let $\left\{ K \right\}_{K \in \Tr}$, the mesh, be a set of closed
polygonal subsets of $\R^d$ with non-empty disjoint interiors such
that $\R^d = \bigcup_{K \in \Tr} K$. The volume ($d$-Lebesgue measure)
of a given cell $K \in \Tr$ is denoted by $|K|$. The supremum of the
diameters of all the cells is denoted by $h$, i.e. $h = \sup_{K \in
\Tr} {\rm diam}(K)$. Two cells $K$ and $L$ are said adjacent if they
aren't disjoint but have disjoint interiors. In this case, we write $K
\sim L$. We assume that, for all pairs $(K,L)$ of adjacent cells, the
intersection $K \cap L$ is included in a hyperplane of $\R^d$. The
surface ($(d-1)$-Lebesgue measure) of the face $K \cap L$ is then
denoted by $|K \cap L|$.

Let $\Delta t > 0$ be the time step of the method. The value
$u_K^n$ intends to approximate the mean value of $u(n \Delta t,\cdot)$
in the cell $K$.
The upwind scheme provides a way to compute such $u_K^n$. 
It is easily obtained by integrating the divergence form of the PDE in
\eqref{trpb}, $\partial_t u + \div (au) - u \div (a) = 0$,
over $[n\Delta t, (n+1)\Delta t] \times K$. We get 
\begin{multline}\label{integ_transp}
\int_K u((n+1)\Delta t,x) dx - \int_K u(n\Delta t,x) dx \\
+ \sum_{L \sim K} \int_{K\cap L} \int_{n\Delta t}^{(n+1)\Delta t}
\left< a(x) , n_{K,L} \right> u(t,x) dt dx - \int_K \int_{n\Delta
t}^{(n+1)\Delta t} u(t,x) \div (a)(x) dt dx = 0, 
\end{multline}
where $n_{K,L}$ is the unit normal vector on $K \cap L$ outward from
$K$. From a numerical point of view, it then seems natural to
compute both an approximate value $u_K^n$ of the mean of
$u(n\Delta t,\cdot)$ in cell the $K$, i.e.
\begin{equation*}
u_K^n \approx |K|^{-1}\int_K u(n\Delta t,x)
dx, 
\end{equation*}
and an approximate value $u_{K,L}^n$ of the mean of $u$ on the
edge $K\cap L$ between the time steps $n$ and $n+1$, i.e.
\begin{equation*}
u_{K,L}^n \approx \Delta t^{-1}|K\cap
L|^{-1} \int_{K\cap L} \int_{n\Delta t}^{(n+1)\Delta t} u(t,x) dt
dx.
\end{equation*}
The quantity $u_{K,L}^n$ is called the {\em numerical flux}. 
Defining $a_{K,L}$ as the mean value of $a$ on the edge $K \cap L$,
i.e. 
\begin{equation*}
a_{K,L} = |K \cap L|^{-1} \int_{K \cap L} a(x) dx, 
\end{equation*}
we get the following approximate version of \eqref{integ_transp}, 
$$
|K| \frac{u_K^{n+1} - u_K^n}{\Delta t} + \sum_{L \sim K}
\left< a_{K,L} , n_{K,L} \right> |K \cap L| u_{K,L}^n - u_K^n \sum_{L
\sim K} \left< a_{K,L} , n_{K,L} \right> |K \cap L| = 0. 
$$
The upwind scheme considers the numerical fluxes $u_{K,L}^n$ as {\em
upwinded}: $u_{K,L}^n = u_K^n$ for $L \in K^+$ and $u_{K,L}^n = u_L^n$
for $L \in K^-$ with 
$$
\begin{array}{l}
K^+ = \left\{ L \sim K, \left< a_{K,L} , n_{K,L} \right> > 0 \right\},
\\
K^- = \left\{ L \sim K, \left< a_{K,L} , n_{K,L} \right> < 0
\right\}. 
\end{array}
$$
This finally gives 
\begin{equation}
\label{eq:scheme}
|K| \frac{u_K^{n+1} - u_K^n}{\Delta t} + \sum_{L \in K^-} \langle
a_{K,L} , n_{K,L} \rangle |K \cap L| \bigl( u_L^n - u_K^n \bigr) = 0,
\ (n,K) \in \N \times \Tr. 
\end{equation}
The numerical initial condition is usually taken
as $u_K^0 =  |K|^{-1} \int_K u^0(x)\, dx $. It is straightforward that
the scheme satisfies the maximum principle under the condition
$$
- \sum_{L \in K^-} \frac{\langle a_{K,L} , n_{K,L} \rangle |K \cap
L|}{|K|} \leq 1, \ K \in \Tr. 
$$
This condition is called the Courant-Friedrichs-Lewy (CFL for short)
condition and is assumed to be satisfied in all the paper. 

\section{Analysis in Dimension 1}\label{d1}

For pedagogical reasons, we first investigate the one-dimensional
framework. As announced in Introduction, the velocity field $a:\R
\rightarrow \R$ is assumed to be bounded and to be $\kappa$-Lipschitz
continuous. In particular, for any starting point $x \in \R$, the
characteristic equation starting from $x$
\begin{equation}
\label{eq:ODEd1}
\partial_t Z (t,x) = - a(Z(t,x)), \ t \geq 0, \quad Z(0,x) = x,
\end{equation}
admits a unique solution. Denoting by $u^0$ the initial condition of
the transport equation, which is assumed to be $\kappa$-Lipschitz
continuous in the whole section, the solution of the transport
equation rewrites
\begin{equation}
\label{eq:reprd1}
u(t,x) = u^0(Z(t,x)), \ (t,x) \in \R_+ \times \R.
\end{equation}

In this section devoted to dimension 1, for every cell $K \in \Tr$,
the volume (length) of $K$ is denoted $\Delta x_K$. The edge
value of $a$ is defined as $a_{K,L} = a(K\cap L)$. The constants
``$C$'' and ``$c$'' used below only depend on $\|a\|_{\infty}$ and
$\kappa$. They are always independent of $\Delta t$, of $h =
\sup_{K\in\Tr} \Delta x_K$, of the time index $n$ and of the random
outcome $\omega$. In particular, the notation $O(x)$, for a given
variable $x$, denotes a quantity bounded by $C x$ for some constant
$C$ only depending on $\|a\|_{\infty}$ and $\kappa$. 

\subsection{Probabilistic Interpretation}

In the one-dimensional framework, the scheme has the form
\begin{equation}
\label{d1:scheme}
\begin{split}
&u^0_K = \frac{1}{\Delta x_K} \int_K u^0(x)dx, \quad K \in \Tr, \\
&u^{n+1}_K = - \sum_{L \in K^-} \frac{a_{K,L} \Delta
t}{\Delta x_K} u^n_{L} + \bigl( 1 + \sum_{L \in K^-}
\frac{a_{K,L} \Delta t}{\Delta x_K} \bigr) u^n_K, \quad n \geq 0,
\quad K \in {\Tr},
\end{split}
\end{equation}
and the following CFL condition is assumed to be in force
\begin{equation}\label{d1:CFL}
- \sum_{L \in K^-} \frac{a_{K,L} \Delta t}{\Delta x_K} \leq 1, \quad K
\in \Tr.
\end{equation}
The geometry of the mesh is simple: each cell $K$ has two
neighbors. When the velocity field $a$ is non-zero in the
cell $K$, there is one and only one cell $L$ in $K^-$ . If $a$ is
positive in $K$, it is the left one; of course, if $a$ is negative, it
is the right one.

We focus for a while on a given cell $K$. By the CFL condition
\eqref{d1:CFL}, all the coefficients
\begin{equation*}
\begin{split}
&p_{K,L} = - \frac{a_{K,L} \Delta t}{\Delta x_K} \ {\rm for } \
L \in K^-, \\
&p_{K,K} = 1 + \sum_{L \in K^-} \frac{a_{K,L}
\Delta t}{\Delta x_K}, \\
&p_{K,L} = 0 \ {\rm for } \ L \in \Tr \setminus \left( K^- \cup K
\right), \\ 
\end{split}
\end{equation*}
are non-negative and may be seen as probability weights.
Henceforth, for a given time step $n \geq 0$, the right-hand side in
\eqref{d1:scheme} may be interpreted as an expectation with respect to
these weights:
\begin{equation*}
u^{n+1}_K = \sum_{L \sim K} p_{K,L} u^n_L.
\end{equation*}

Intuitively, this means that we are choosing one cell among
$K \cup \{L\sim K\}$, $K$ being fixed, with the probability weights
$p_{K,K}$ and $p_{K,L}$ for $L\sim K$. To make this idea rigorous,
we introduce a probability space $(\Omega,{\mathcal
A},{\mathbb P})$ as well as a random variable $\xi : \Omega
\rightarrow K \cup \{L\sim K\}$ such that ${\mathbb P}\{\xi=L\} =
p_{K,L}$ for any $L \sim K$ and ${\mathbb P}\{\xi=K\} =
p_{K,K}$. Then, the $(n+1)^{\scriptsize{th}}$ step of the numerical
scheme on the cell $K$ can be written in the following way:
\begin{equation}
\label{eq:1} u^{n+1}_K = \sum_{L \sim K} p_{K,L} u^n_L = {\mathbb E}
\bigl[ u^n_{\xi} \bigr].
\end{equation}
This relationship provides a probabilistic interpretation for the
one step dynamics of the numerical scheme. We are to iterate this
procedure.

The probabilistic dynamics between times $n$ and $n+1$ just depend on
the \emph{starting cell} $K$. In the theory of stochastic processes, this
property is typical of \emph{Markovian dynamics}. Indeed, the family of
probability weights $(p_{K,L})_{K,L \in {\Tr}}$ defines a
stochastic matrix of infinite dimension (all the entries of the matrix
are non-negative and the sums of the entries of a same line are all
equal to 1). This stochastic matrix corresponds to the transition
matrix of a Markov chain. Up to a modification of the underlying
probability space, there exists a sequence $(K_n)_{n \geq 0}$ of
random variables taking values into the set of cells as well as a
collection of probability measures $({\mathbb P}_{K})_{K \in \Tr}$,
indexed by the cells, such that, under each ${\mathbb P}_K$, $(K_n)_{n
\geq 0}$ is a Markov chain with rates $(p_{K,L})_{K,L \in {\Tr}}$
starting from $K_0=K$. In other words,
\begin{equation*}
\forall n \geq 0, \ {\mathbb P}_{K}\{K_{n+1}=L|K_n=K\} = p_{K,L}, \
{\mathbb P}_{K}\{K_0=K\} =1.
\end{equation*}
The behavior of the chain $(K_n)_{n \geq 0}$ is as follows: if the
velocity is positive in the cell $K_n$, then the probability
$p_{K_n,L}$ vanishes if $L$ is the right neighbor of $K_n$, so that
the chain can either stay in $K_n$ or jump to the left.

Now, we can interpret \eqref{eq:1} in a different way:
\begin{equation*}
u^{n+1}_K = {\mathbb E}_K \bigl[ u^n_{K_1} \bigr],
\end{equation*}
where ${\mathbb E}_K$ denotes the expectation associated with
${\mathbb P}_K$. This means that $u^{n+1}_K$ is the expectation of
$u^n$ in the random cell $K_1$ occupied by the Markov chain, which
started one time step before in $K$.
We can also write for any integer $i \geq 0$
\begin{equation*}
u^{n+1}_{K_i} = {\mathbb E}_K \bigl[ u^n_{K_{i+1}} | K_0,\dots,
K_i\bigr] \quad {\mathbb P}_K {\rm -almost \ surely.} 
\end{equation*}
(When conditioning with respect to $K_0,\dots,K_i$, the past before
$i-1$ doesn't play any role, and the chain restarts, afresh, at time
$i$ from $K_i$.) In what follows, we denote the conditional
expectation ${\mathbb E}_K[\cdot | K_0,\dots,K_n]$ by ${\mathbb
E}^n_K[ \cdot]$. We also omit to specify that such a conditional
expectation is computed under ${\mathbb P}_K$. Now, we are able to
iterate 
the procedure in \eqref{eq:1}:
\begin{equation*}
\begin{split}
u^{n+1}_K = {\mathbb E}_K \bigl[ u^n_{K_1} \bigr] &= {\mathbb E}_K
\bigl[ {\mathbb E}^1_K \bigl[ u^{n-1}_{K_2} \bigr]] \\
&= \cdots = {\mathbb E}_K \bigl[ {\mathbb E}^1_K \bigl[ \cdots
{\mathbb E}^{n}_K \bigl[u^0_{K_{n+1}} \bigr] \bigr] \bigr] = {\mathbb
E}_K \bigl[u^0_{K_{n+1}} \bigr].
\end{split}
\end{equation*}
We have proved the following representation for the numerical
solution $u^n$:
\begin{theorem}\label{th:feynman_kac}
Under the above notations, the numerical
solution $u^n_K$ at time $n$ and in the cell $K$ has the form:
\begin{equation*}
u^n_K = {\mathbb E}_K \bigl[ u^0_{K_n} \bigr].
\end{equation*}
\end{theorem}
The representation given by Theorem \ref{th:feynman_kac} is a
backward Kolmogorov formula for the numerical scheme. Generally
speaking, the backward Kolmogorov formula provides a representation
for the solution of the heat equation in terms of the mean value of
the initial condition computed with respect to the paths of the
Brownian motion (or of a diffusion process). In this framework, the
paths of the Brownian motion appear as random characteristics. In our
own setting, the Markov chain $(K_n)_{n \geq 0}$ almost plays the same
role.

We say ``almost plays'' because the sequence $(K_n)_{n \geq 0}$ is
not a sequence of points as the Brownian motion is. Actually, we
have to associate with each random cell $K_n$ a random point $X_n$
($X_n$ being ideally in $K_n$) to obtain a random characteristic
$(X_n)_{n \geq 0}$.

The choice of these points is crucial. In what follows, we choose
$X_n$ as the entering point in the cell $K_n$. This means that
\begin{equation*}
X_n = X_{n-1} \ {\rm if} \ K_n = K_{n-1}, \ X_n = K_n \cap K_{n-1} \
{\rm if} \ K_n \not = K_{n-1}.
\end{equation*}
The above definition holds for $n \geq 1$. The position of the
initial point $X_0$ inside $K_0$ has to be specified. If $a(x)>0$ for
all $x \in K_0$, we choose $X_0$ as the right boundary of
$K_0$. (Indeed, the right boundary plays in this case the role of the
entering point since the velocity is positive.) If $a(x)<0$ for all $x
\in K_0$, we choose $X_0$ as the left boundary. If $\exists x \in K_0$
such that $a(x)=0$, we choose $X_0$ as the middle of $K_0$.

What is important is that the sequence $(X_n)_{n \geq 0}$ is adapted
to the filtration generated by $(K_n)_{n \geq 0}$, i.e.
$(\sigma(K_0,\dots,K_n))_{n \geq 0}$: knowing the paths
$(K_0,\dots,K_n)$, one knows the positions of the points
$(X_0,\dots,X_n)$.

The reader may wonder about this specific choice for the sequence
$(X_n)_{n \geq 0}$. Assume that the velocity $a$ is non-zero, say
for example positive, in the cell $K_n$. By continuity, it is
positive in the neighborhood of $K_n$: the chain goes from the right
to the left in this area of the space. As a by-product, the entering
point in the cell $K_n$ is the right boundary of $K_n$. In this
case, $X_{n+1}$ is either $X_n$ itself or the left boundary of
$K_{n}$, which is the right boundary of $K_n^-$ (this set of cells is
in the present case a singleton and we identify it with its element,
as well as for $K_n^+$ in the following when the velocity is away from
0), so that
\begin{equation*}
{\mathbb E}^n_K \bigl[X_{n+1} - X_n \bigr] = - \Delta x_{K_n}
p_{K_n,K_n^-} = - a_{K_n,K_n^-} \Delta t.
\end{equation*}
(Indeed, the probability that $X_{n+1}$ is the right boundary of
$K_n^-$ is given by the probability of jumping from $K_n$ to $K_n^-$.)
In other words, the mean displacement from $X_n$ to $X_{n+1}$, knowing
the past, is exactly driven by the velocity field $-a$. This is very
important: loosely speaking, the speed of the random characteristic is
given by the velocity field of the characteristic equation itself!
Here is a more precise statement. 
\begin{prop}\label{prop:mean_speed}
For every $ n\geq 0$,
\begin{equation*}
{\mathbb E}^n_K \bigl[X_{n+1} - X_n \bigr] = -a (X_n) \Delta t + O(h
\Delta t).
\end{equation*}
\end{prop}
\noindent {\bf Proof.} If $a(x)>0$ for all $x \in K_n$, then $X_n$ has
to be the right boundary of $K_n$. (By definition of $X_0$, this is
true until the first jump of the chain. After the first jump, this is
still true since the chain cannot come from the left by positivity of
$a$.) Moreover, starting from $K_n$, the chain cannot move to the
right since $a_{K_n,K_n^+}=a(X_n) >0$. Hence, $X_{n+1}$ has to be
either $X_n$ or the left boundary of $K_{n}$. As done above,
\begin{equation*}
{\mathbb E}^n_K \bigl[X_{n+1} - X_n \bigr] = - \Delta x_{K_n}
p_{K_n,K_n^-} = - a_{K_n,K_n^-} \Delta t = - a(X_n) \Delta t + O(h
\Delta t)
\end{equation*}
by the Lipschitz property of $a$. The same argument holds when $a(x)
<0$ for all $x \in K_n$.

If $a(x)=0$ for some $x \in K_n$, then $a(X_n) = O(h \Delta t)$ by the
Lipschitz property of $a$. Moreover, the probability of moving to the
right is equal to $\max(a_{K_n,K_n^+},0) \times \Delta t/ \Delta
x_{K_n} = O(\Delta t)$. The same holds for the probability of moving
to the left. When moving, the displacement is bounded by $h$ so that
the result is still true.
\qed
\vspace{5pt}

\noindent
{\bf Remark.}
The necessity of choosing the point $X_n$ as the entering point in the
cell $K_n$ is related to the well-known fact that the upwind
scheme is consistent (in the finite difference sense) with the
transport equation provided that the control points for every cell are
chosen on the right if the velocity is positive and, conversely, on
the left if the velocity is negative: see \cite{EymGalHer00}.

\subsection{A First Example: $a$ and $\Delta x$ constant}

To explain our strategy, we first focus on the very simple case
where both $a$ and $\Delta x$ are constant: $a(x)=a$ for all $x \in
{\mathbb R}$ and $\Delta x_K = h$ for all $K \in \Tr$. Without
loss of generality, we can assume that $a$ is positive so that the
random characteristic goes from the right to the left. In this
setting, the transition probabilities are of the form
\begin{equation*}
\forall K \in \Tr, \ p_{K,K^-} = \frac{a \Delta t}{h}, \ p_{K,K} = 1-
\frac{a \Delta t}{h}.
\end{equation*}
The probability of jumping from one cell to another doesn't depend
on the current state of the random walk. From a probabilistic point
of view, this amounts to say that the sequence $(X_{n+1}-X_n)_{n\geq
0}$ (with $X_0$ equal to the right boundary of the initial cell) is a
sequence of Independent and Identically Distributed (IID in short)
random variables under ${\mathbb P}_K$, whatever $K$ is. The common
law of these variables is given by
\begin{equation*}
{\mathbb P}_K\{ X_{n+1} - X_n = - h\} = 1 -{\mathbb P}
\{X_{n+1} - X_n = 0\} = \frac{a \Delta t}{h}.
\end{equation*}
In particular, we recover a stronger version of Proposition
\ref{prop:mean_speed} (``stronger'' means that there is no $O(h \Delta
t)$): 
\begin{equation*}
{\mathbb E}_K \bigl[X_{n+1} - X_n \bigr] = - a \Delta t.
\end{equation*}
In particular, the mean trend of the random characteristic is exactly
driven by the velocity $-a$, that is by the
mapping $t \mapsto X_0 - a t$, which corresponds to the
characteristic of the transport equation with $X_0$ as initial
condition, i.e. $Z(t,X_0)$ (see \eqref{eq:ODEd1}). At this stage, we
understand that the order of the numerical scheme is deeply related to
the fluctuations of the random characteristic around its mean trend,
that is around the deterministic characteristic. Indeed, for any
starting cell $K$, we have
\begin{equation*}
u_{K}^n = {\mathbb E}_K \bigl[u^0_{K_n} \bigr] = {\mathbb E}_K
\bigl[u^0(X_n) \bigr] + O(h),
\end{equation*}
where $O(h)$ only depends on the Lipschitz constant of the initial
condition $u^0$ and is independent of the initial cell $K$. Thus 
\begin{equation}
\label{d1:100}
u_K^n = {\mathbb E}_K \bigl[ u^0 \bigl( X_0 - a n \Delta t +
\sum_{k=0}^{n-1} (X_{k+1} - X_k + a \Delta t) \bigr) \bigr] + O(h).
\end{equation}
By \eqref{eq:reprd1}, for all $x \in K$, 
\begin{equation*}
\begin{split}
u_K^n - u(n\Delta t,x)
&= {\mathbb E}_K \bigl[ u^0 \bigl( X_0 - a n \Delta t +
\sum_{k=0}^{n-1} (X_{k+1} - X_k + a \Delta t) \bigr) \bigr]  -
u(n\Delta t,x) + O(h) \\ 
&= {\mathbb E}_K \bigl[ u^0 \bigl( X_0 - a n \Delta t +
\sum_{k=0}^{n-1} (X_{k+1} - X_k + a \Delta t) \bigr)   -
u(n\Delta t,X_0) \bigr] + O(h)\\ 
&= {\mathbb E}_K \bigl[ u^0 \bigl(X_0 - a n \Delta t +
\sum_{k=0}^{n-1} (X_{k+1} - X_k + a \Delta t) \bigr) - u^0(X_0 - a n
\Delta t) \bigr] + O(h). 
\end{split}
\end{equation*}
Using again the Lipschitz continuity of $u^0$,
we deduce by the Cauchy-Schwarz inequality that, 
\begin{equation*}
\begin{split}
|u_K^n - u(n\Delta t,x)| &\leq \kappa {\mathbb E}_K \bigl[
\bigl| \sum_{k=0}^{n-1} (X_{k+1} - X_k + a \Delta t) \bigr| \bigr]
+ O(h) \\
&\leq\kappa {\mathbb E}_K \bigl[ \bigl| \sum_{k=0}^{n-1} (X_{k+1} -
X_k + a \Delta t) \bigr|^2 \bigr]^{1/2} + O(h) 
\end{split}
\end{equation*}
for all $x \in K$. 
The last expectation is nothing but the variance of the sum of the
random variables $(X_{k+1}-X_k)_{0 \leq k \leq n-1}$ under ${\mathbb
P}_K$, i.e.
\begin{equation*}
{\mathbb E}_K \bigl[ \bigl| \sum_{k=0}^{n-1} (X_{k+1} -
X_k + a \Delta t) \bigr|^2 \bigr] = {\mathbb V}_K \bigl[
\sum_{k=0}^{n-1} (X_{k+1} - X_k) \bigr].
\end{equation*}
It is well-known that the variance of the sum of independent random
variables is equal to the sum of the variances of the variables. We
deduce that, for all $x \in K$, 
\begin{equation*}
|u_K^n - u(n\Delta t,x)| \leq \kappa \bigl[ n {\mathbb
V}_K(X_1-X_0) \bigr]^{1/2} + O(h).
\end{equation*}
The common variance is equal to
\begin{equation*}
{\mathbb V}_K(X_1 - X_0) = h^2 p_{K,K^-} - a^2 \Delta t^2 =
a\Delta t \bigl( h - a \Delta t \bigr).
\end{equation*}
Note that the CFL condition guarantees that the right-hand side above
is non-negative, which ensures that the equality is meaningful.
We thus recover a well-known estimate for the $L^{\infty}$-error
induced by the upwind scheme:
\begin{prop}
\label{prop:estimate_d1}
Assume that $a(x) = a$ and $\Delta x_K$ are constant and that $u^0$ is
bounded and $\kappa$-Lipschitz continuous. Then, at any time $n \geq
0$, 
\begin{equation*}
\sup_{K \in {\Tr}} \left|\left| u_K^n - u(n \Delta
t,\cdot)\right|\right|_{L^\infty(K)} \leq \kappa \bigl( n a \Delta
t ( h - a \Delta t) \bigr)^{1/2} + O( h). 
\end{equation*}
\end{prop}

\subsection{General Case: $a(x)$ and $\Delta x_K$ non constant}

Our strategy is sharp enough to obtain the analogue of Proposition
\ref{prop:estimate_d1} when $a$ does depend on $x$ and $\Delta x_K$ on
$K$. The main difference here is that $X_0$ can be either the right
boundary of $K_0$ or the left boundary or the
barycenter of $K_0$, according to the sign of the velocity in $K_0$. 
Following the previous subsection, the difference between the
numerical and the true solutions at time $n \geq 0$ on a cell $K$ is
given by 
\begin{equation*}
u^n_K - u(n \Delta t,x) = {\mathbb E}_K \bigl[ u^0_{K_n} \bigr] -
u^0(Z(n \Delta t, x)) = {\mathbb E}_K \bigl[ u^0_{K_n} -
u^0(Z(n \Delta t, x)) \bigr]
\end{equation*}
for all $x \in K$. 
As above, the Lipschitz property yields
\begin{equation*}
u^0_{K_n} = u^0(X_n) + O (h).
\end{equation*}
Again, the term $O(h)$ is uniform with respect to the starting cell
$K$, to the time index $n$, to the parameters $\Delta t$ and $h$ and
to the underlying draw $\omega \in \Omega$. By Gronwall's lemma, we
control the distance between $Z(n \Delta t,x)$ and $Z(n \Delta
t,X_0)$, so that 
\begin{equation}
\begin{split}
\label{eq:3} |u^n_K - u(n \Delta t,x)| &\leq \kappa {\mathbb
E}_K \bigl[ |X_n - Z(n \Delta t, x)| \bigr] + O (h) 
\\
&= \kappa {\mathbb
E}_K \bigl[ |X_n - Z(n \Delta t, X_0)| \bigr] + O(h) \exp(\kappa
n\Delta t). 
\end{split}
\end{equation}
We have
\begin{equation}
\label{d1:200}
X_n - Z(n \Delta t, X_0) = X_n - X_0 + \sum_{k=0}^{n-1} \int_{k
\Delta t}^{(k+1) \Delta t} a(Z(s, X_0))ds.
\end{equation}
By Proposition \ref{prop:mean_speed},
\begin{equation}
\label{d1:201}
X_n - X_0 = \sum_{k=0}^{n-1} \bigl( X_{k+1} - X_k \bigr) = - \Delta t
\sum_{k=0}^{n-1} a(X_k) + M_n + O(n h \Delta t ),
\end{equation}
with
\begin{equation}
\label{d1:205}
M_n = \sum_{k=0}^{n-1} \bigl(X_{k+1} - X_k - {\mathbb
E}_K^k(X_{k+1}-X_k) \bigr) \quad (M_0=0).
\end{equation}
By the boundedness and the Lipschitz continuity of $a$,
\begin{equation}
\label{d1:202}
\sum_{k=0}^{n-1} \int_{k \Delta t}^{(k+1) \Delta t} a(Z(s, X_0))ds =
\Delta t \sum_{k=0}^{n-1} a\bigl( Z(k \Delta t, X_0) \bigr) + O(n
\Delta t^2).
\end{equation}
Plugging \eqref{d1:201} and \eqref{d1:202} into \eqref{d1:200}, we
obtain
\begin{equation*}
|X_n - Z(n \Delta t, X_0)| \leq \kappa \Delta t \sum_{k=0}^{n-1}
|X_{k} - Z(k \Delta t, X_0)| + |M_n| + O ( n h \Delta t + n \Delta
t^2 ).
\end{equation*}
Taking the expectation of each term and applying Gronwall's lemma,
\begin{equation}
\label{d1:300}
{\mathbb E}_K \bigl[ |X_n - Z(n \Delta t, X_0)| \bigr] \leq \bigl[
{\mathbb E}_K[|M_n|] + O( n h \Delta t + n \Delta t^2 )\bigr] \exp(
\kappa n \Delta t).
\end{equation}

As in the case where both the velocity and the spatial step are
constant, the process $(M_n)_{n \geq 0}$ represents the fluctuations
of the random characteristic around a discretized version of the
deterministic characteristic. (See \eqref{d1:201}.) In the
probabilistic theory, it is a martingale on $(\Omega,{\mathcal
A},{\mathbb P}_K)$, i.e., at any time $n \geq 0$,
$M_n$ is $\sigma(K_0,\dots,K_n)$-measurable and
${\mathbb E}_K^n[M_{n+1}] = M_n$. This property just follows from
\eqref{d1:205}. 

Since $M_0=0$, the expectation of $M_n$, for $n \geq 1$, is given by
${\mathbb E}_K[M_n]={\mathbb E}_K[{\mathbb E}_K^{n-1}(M_n)] = {\mathbb
E}_K[M_{n-1}] = \dots = M_0=0$. The mean trend of a martingale
starting from zero is null. To estimate the fluctuations, we compute
the second order moment. Setting $\Delta M_j = M_{j+1} - M_{j}$ for
all $j \geq 0$, the martingale property yields ${\mathbb
E}_K^{j}[\Delta M_j]=0$, so that, for all $n \geq 1$,
\begin{equation}
\label{eq:4}
\begin{split}
{\mathbb E}_K (M_n^2) &= \sum_{k=0}^{n-1} {\mathbb E}_K \bigl[\Delta
M_k^2 \bigr] + 2 \sum_{0 \leq i < j \leq n-1} {\mathbb E}_K \bigl[
\Delta M_i \Delta M_j \bigr] \\
&= \sum_{k=0}^{n-1} {\mathbb E}_K \bigl[\Delta
M_k^2 \bigr] + 2 \sum_{0 \leq i < j \leq n-1} {\mathbb E}_K \bigl[
 {\mathbb E}_K^j (\Delta M_i \Delta M_j) \bigr] \\
&= \sum_{k=0}^{n-1} {\mathbb E}_K \bigl[\Delta M_k^2 \bigr] + 2
\sum_{0 \leq i < j \leq n-1} {\mathbb E}_K \bigl[ \Delta M_i
{\mathbb E}_K^j (\Delta M_j) \bigr] \\
&= \sum_{k=0}^{n-1} {\mathbb E}_K \bigl[\Delta M_k^2 \bigr].
\end{split}
\end{equation}
It remains to compute the expectation of the increments $(\Delta
M_k^2)_{k \geq 1}$.

We first prove that ${\mathbb E}_K^k[\Delta M_k^2] = O (h \Delta
t)$. Since $\Delta M_k = X_{k+1} - X_k - {\mathbb
E}_K^k(X_{k+1}-X_k)$, we have (${\mathbb V}_K^k$ denotes the
conditional variance knowing $K_0,\dots,K_k$ under ${\mathbb P}_K$)
\begin{equation}
\label{eq:varianced1} 
\begin{split}
{\mathbb E}_K^k [\Delta M_k^2] &= {\mathbb V}_K^k [X_{k+1} - X_k] = 
{\mathbb E}_K^k [(X_{k+1} - X_k)^2] - \bigl({\mathbb E}_K^k [X_{k+1} -
X_k] \bigr)^2 
\\
&\leq {\mathbb E}_K^k \bigl[ (X_{k+1} - X_k)^2 \bigr].
\end{split}
\end{equation}
Knowing the position of the chain at time step $k$,
the conditional probability of jumping is bounded by
$(\max(a_{K_k,K_k^-},0)+\max(a_{K_k,K_k^+},0))\Delta t/\Delta
x_{K_n}$. When jumping, the distance between $X_k$ and $X_{k+1}$ is
always bounded by $\Delta x_{K_k}$. Hence, ${\mathbb E}_K^k [ (X_{k+1}
- X_k)^2 ] = O(h \Delta t)$. Taking the expectation, we deduce that
${\mathbb E}_K [\Delta M_k^2] = O(h \Delta t)$.

By \eqref{eq:4}, we deduce
\begin{equation*}
{\mathbb E}_K(M_n^2) = O(n h \Delta t).
\end{equation*}
By \eqref{eq:3} and \eqref{d1:300} and by the Cauchy-Schwarz
inequality, we deduce

\begin{prop}\label{prop:d1_general}
Under the assumptions introduced in the beginning of Section \ref{d1},
there exists a constant $C \geq 0$, such that at any time $n \geq 0$,
\begin{equation*}
\sup_{K \in {\Tr}} \left|\left| u^n_K - u(n \Delta t,\cdot)
\right|\right|_{L^\infty(K)} \leq C \bigl( (n h \Delta t)^{1/2} + n h
\Delta t+ n \Delta t^2 + h \bigr) \exp(\kappa n \Delta t).
\end{equation*}
\end{prop}

\subsection{Interpretation by the Central Limit Theorem}

This section only concerns the special case with constant velocity on
a uniform mesh. It provides a finer result in this simplified case, by
the use of the central limit theorem. This analysis will not be
performed in higher dimension. We again assume that $a(x) = a > 0$ and
$\Delta x_K$ is constant. We also reinforce the CFL condition, asking
$a \Delta t < h$. (This is not a restriction. When $a \Delta t =h$,
the term of order $1/2$ vanishes in Proposition \ref{prop:estimate_d1}
and the error is of order $1$: this case is trivial.) As explained
above, the order of the numerical scheme is given by the order of the
fluctuations of the random characteristic around its mean trend. In
the specific setting where both $a$ and $\Delta x$ are constant, the
random characteristic corresponds to a random walk with IID
increments: by the elementary theory of stochastic processes, we know
that the fluctuations of the walk around its mean trend are governed
by the Central Limit Theorem (CLT in short). (See \cite[Chapter III,
\S 3]{shiryaev} for the standard version of the CLT and \cite[Chapter
2, Theorem 4.20]{karatzas:shreve} for the functional version in the
case of a simple random walk.) We deduce that the fluctuations are of
diffusive type, that is they correspond, asymptotically, to the
fluctuations of a Brownian motion around the origin. This means that
the random characteristics can be seen, asymptotically, as the paths
of a Brownian motion, with a non-standard variance. (That is, the
variance of the Brownian motion at time $t$ isn't equal to $t$, but to
a constant times $t$, the constant being independent of $t$. In our
framework, the constant is proportional to $h$.) From an analytical
point of view, we are saying that the numerical solution is very close
to the solution of a second order parabolic equation: this is nothing
but the numerical diffusive effect.

We specify this idea. The variables $(X_{n+1}-X_n + a \Delta t)_{n
\geq 0}$ are IID, with zero as mean and $a \Delta t (h - a
\Delta t)$ as variance. By the CLT, we know that
\begin{equation}\label{d1:101}
\bigl( n a \Delta t (h - a \Delta t) \bigr)^{-1/2} \sum_{k=0}^{n-1}
(X_{k+1} - X_k + a \Delta t) \, \Rightarrow \, {\mathcal N}(0,1) \quad
{\rm as }\ n \rightarrow + \infty,
\end{equation}
on each $(\Omega,{\mathcal A},{\mathbb P}_K)$, $K \in {\Tr}$.
The notation $\Rightarrow$ stands for the convergence in
distribution. (In short, for a family of random variables $(W_n)_{n \in
\{0,\dots,+\infty\}}$ on $(\Omega,{\mathcal A},{\mathbb P})$, we say
that $W_n \Rightarrow W_{\infty}
$ if ${\mathbb E}[\varphi(W_n)] \rightarrow {\mathbb
E}[\varphi(W_{\infty})]$ for any bounded continuous function
$\varphi$.) The notation ${\mathcal N}(0,1)$ stands for the reduced
centered Gaussian law.

Plugging \eqref{d1:101} into \eqref{d1:100}, we deduce that
\begin{equation*}
u^n_K \approx {\mathbb E}_K \bigl[ u^0
\bigl(X_0 - an \Delta t+ (n a \Delta t (h- a \Delta t))^{1/2} W
\bigr) \bigr] + O(h)
\end{equation*}
for $n$ large. Above, $W$ denotes a reduced centered Gaussian random
variable. The symbol $\approx$ means that both sides are close. This
point will be specified below.

Expliciting the density of the Gaussian law, we can write
\begin{equation}
\label{d1:103}
\begin{split}
&{\mathbb E}_K \bigl[ u^0
\bigl(X_0 - an \Delta t+ (n a \Delta t (h- a \Delta t))^{1/2} W
\bigr) \bigr] \\
&= \int_{\mathbb R} u^0(x_{K,K^+}-y) \exp \bigl[ - \frac{
\bigl(y - a n \Delta t)^2}{2 a (h- a \Delta t) n \Delta t} \bigr]
\frac{dy} {[2 \pi a (h - a \Delta t) n \Delta t]^{1/2}},
\end{split}
\end{equation}
where $x_{K,K^+}$ stands for the unique point in the intersection 
of $K$ and $K^+$, i.e. the right boundary of the cell $K$ when $a>0$.
According to \cite[Chapter 1]{friedman} (with
$a=0$, the generalization to $a \not =0$ being trivial), we recognize
the value at time $n \Delta t$ and at point $x_{K,K^+}$ of the
solution $v$ to the Cauchy problem: 
\begin{equation}
\label{d1:104}
\partial_t v + a \partial_x v - \frac{a ( h - a \Delta t)}{2}
\partial^2_{x,x} v = 0, \ t \geq 0, \ x \in \R,
\end{equation}
with $u^0$ as initial condition. Finally, we can say that $u^n_K$ is
close to $v(n \Delta t,\cdot)$ in the cell $K$.

Of course, we have to say what ``close'' means! This question is
related to the rapidity of convergence in the CLT, that is the
rapidity of convergence in \eqref{d1:101}. The main result in this
direction is the Berry-Esseen Theorem. (See \cite[Chapter III, \S
11]{shiryaev}.) In what follows, we use a refined version of it. (See
\cite[Chapter V, \S 4, Theorem 14]{petrov}.) Denoting by $F_n$ the
cumulative distribution function
\begin{equation*}
\forall z \in \R, \ F_n(z) = {\mathbb P}_K \bigl\{ \bigl(n a \Delta t
(h- a \Delta t) \bigr)^{-1/2} \sum_{k=0}^{n-1} (X_{k+1} - X_k + a
\Delta t) \leq z \bigr\}, 
\end{equation*}
and by $\Phi$ the cumulative distribution function of the ${\mathcal
N}(0,1)$ law, we have, for all $z \in \R$,
\begin{equation*}
|F_n(z) - \Phi(z)| \leq C n^{-1/2} \bigl(
a \Delta t(h- a \Delta t) \bigr)^{-3/2} {\mathbb E}_K \bigl[ |X_1-X_0+
a \Delta t|^3 \bigr] (1+|z|)^{-3},
\end{equation*}
for some universal constant $C>0$.
The moment of order three is given by
\begin{equation*}
\begin{split}
{\mathbb E}_K \bigl[ |X_1-X_0 + a \Delta t|^3 \bigr]
&= \frac{a \Delta t}{h} ( h - a \Delta t)^3 + \bigl(1 - \frac{a
\Delta t}{h} \bigr) (a \Delta t)^3 \\
&= a \Delta t \bigl( 1 - \frac{a \Delta t}{h} \bigr) \bigl[ ( h - a
\Delta t)^2 + (a \Delta t)^2 \bigr].
\end{split}
\end{equation*}
Hence, for all $z \in \R$,
\begin{equation}
\label{d1:102}
|F_n(z) - \Phi(z)| \leq C \bigl(n
a \Delta t (h- a \Delta t) \bigr)^{-1/2} h^{-1}
\bigl[ ( h - a \Delta t)^2
+ (a \Delta t)^2 \bigr] (1+|z|)^{-3}.
\end{equation}
By \eqref{d1:100},
\begin{equation*}
u^n_K = \int_{\R} u^0\bigl(x_{K,K^+} - an \Delta t
+ (n a \Delta t (h- a \Delta t))^{1/2} y \bigr) dF_n(y) + O(h),
\end{equation*}
where the integral in the right-hand side is a
Lebesgue-Stieltjes integral. (See \cite[Chapter II, \S 6]{shiryaev}.)

Assume for a while that the support of $u^0$ is compact.
Performing an integration by parts (see \cite[Chapter II, \S 6,
Theorem 11]{shiryaev}), we obtain (since $u^0$ is Lipschitz
continuous) 
\begin{equation*}
u^n_K = (n a \Delta t (h- a \Delta t))^{1/2} \int_{\R}
\frac{du^0}{dy}\bigl(x_{K,K^+} - an \Delta t + (n a \Delta t (h- a
\Delta t))^{1/2} y \bigr) F_n(y) dy + O(h),
\end{equation*}
Plugging \eqref{d1:102} in this equality, we obtain
\begin{equation*}
\begin{split}
&\bigl| u^n_K - (n a \Delta t (h- a \Delta t))^{1/2} \int_{\R}
\frac{du^0}{dy} \bigl(x_{K,K^+} - an \Delta t + (n a \Delta t (h- a
\Delta t))^{1/2} y \bigr) \Phi(y) dy \bigr| \\
& \hspace{15pt} \leq C \kappa h^{-1} \bigl[( h - a \Delta t)^2
+ (a \Delta t)^2 \bigr] \int_{\R} (1+|y|)^{-3} dy + O \bigl(h).
\end{split}
\end{equation*}
Performing a new integration by parts and then a change of
variable, we see that the left-hand side is equal to $|u_K^n - v(n
\Delta t,x_{K,K^+})|$, where $v$ is the solution of
the Cauchy problem \eqref{d1:104}. Using a standard troncature
argument, we can easily get rid of the assumption made on the support
of $u^0$. We finally claim
\begin{prop}\label{prop:TLC}
Assume that $a(x)$ and $\Delta x_K$ are constant and that $u^0$ is
$\kappa$-Lipschitz continuous. Then, there exists a constant $C>0$
such that, at time any time $n \geq 0$, 
\begin{equation*}
\sup_{ K \in {\Tr}} \left|\left| u^n_K - v(n \Delta t,\cdot)
\right|\right|_{L^\infty(K)} \leq C h,
\end{equation*}
where $v$ stands for the solution of the Cauchy problem
\eqref{d1:104} with $u^0$ as initial condition.
\end{prop}
It is remarkable that the bound doesn't depend on
$(n,\Delta t)$. The result is still true for $a \Delta t= h$. (See
Proposition \ref{prop:estimate_d1}.)

\section{Principle of the Analysis in higher dimension. Application to
a Simple Case}\label{sec:2}

We here present the basic ingredients for the analysis
in dimension $d$ greater than two.

As above, the velocity field $a$ is assumed to be bounded and
$\kappa$-Lipschitz continuous. The regularity of the initial condition
$u^0$ will be specified below.

The characteristics of the transport equation are still denoted by
$(Z(t,x))_{t \geq 0}$, $x \in \R^d$, see Equation \eqref{carac}.

As in \cite{merlet:vovelle}, the cells are assumed to be uniformly non
flat, i.e. they satisfy, in a strong sense, the converse of the
isoperimetric inequality:
\begin{equation}
\label{ass:iso} \exists \alpha >0, \ \forall K \in {\Tr}, \
\sum_{L \sim K} |K \cap L| \leq \alpha |K| h^{-1}.
\end{equation}
By the standard isoperimetric inequality, this is equivalent to the
existence of $\beta >1$ such that $|K| \geq \beta^{-1} h^d$ and
$\sum_{L \sim K} |K \cap L| \leq \beta h^{d-1}$ for all $K \in
{\Tr}$.

The mean velocity on a cell $K$ is denoted by
\begin{equation*}
a_K = |K|^{-1} \int_{K} a(x) dx.
\end{equation*}

The constants ``$C$'' and ``$c$'' below may depend on
$\|a\|_{\infty}$, $\alpha$, $\beta$, $\kappa$ and $d$. As in
dimension one, they are always independent of $\Delta t$, of $h$, of
the current time index $n$ and of the random outcome $\omega$.

\subsection{Stochastic Representation of the Scheme}

The expression of $u_K^{n+1}$ given by the upwind Scheme is (see
\eqref{eq:scheme}) 
\begin{equation*}
u_K^{n+1} = - \sum_{L \in \Km} \frac{\langle a_{K,L},n_{K,L}\rangle
\Delta t |K \cap L|}{|K|} u^n_L + \left(1 + \sum_{L \in
\Km}\frac{\langle a_{K,L},n_{K,L}\rangle \Delta t |K \cap L|}{|K|}
\right) u^n_K,
\end{equation*}
where $\Km = \{L \sim K, \ \langle a_{K,L},n_{K,L} \rangle <0\}$,
$u^0_K$ being given by
\begin{equation*}
u^0_K = |K|^{-1} \int_K u^0(x) dx.
\end{equation*}
In this framework, the CFL condition has the form
 \begin{equation*}
\forall K \in {\Tr}, \
- \sum_{L \in \Km}\frac{\langle
a_{K,L},n_{K,L}\rangle \Delta t |K \cap L|}{|K|} \leq 1
\end{equation*}
and is assumed to be satisfied in all the following.

As in dimension one, the coefficients
\begin{equation*}
\begin{split}
&p_{K,L} = - \frac{\langle
a_{K,L},n_{K,L}\rangle \Delta t |K \cap L|}{|K|} \ {\rm for} \ L \in
\Km, \\
&p_{K,K} = 1 + \sum_{L \in \Km}\frac{\langle
a_{K,L},n_{K,L}\rangle \Delta t |K \cap L|}{|K|}, \\
& p_{K,L} = 0 \ {\rm for} \ L \in \Tr \setminus \left(K^- \cup
K\right), 
\end{split}
\end{equation*}
can be interpreted as the probability transitions of a Markov chain
with values in the set of cells. Again, we can find a measurable
space $(\Omega,{\mathcal A})$, a sequence $(K_n)_{n \geq 0}$ of
measurable mappings from $(\Omega,{\mathcal A})$ into ${\Tr}$
as well as a family $({\mathbb P}_K)_{K \in {\Tr}}$ of
probability measures on $(\Omega,{\mathcal A})$, indexed by the cells,
such that, for every cell $K \in {\Tr}$, $(K_n)_{n \geq 0}$ is
a Markov chain with $K$ as initial condition and $(p_{K,L})_{K,L \in
{\Tr}}$ as transition probabilities.
As in dimension one, the chain $(K_n)_{n \geq 0}$ goes against the
velocity field $a$: for $n \geq 0$, either $K_{n+1}$ is equal to $K_n$
or $K_{n+1}$ belongs to $ K_n^-$. Similarly, for $n \geq 1$, either
$K_{n-1}$ is equal to $K_n$ or $K_{n-1}$ belongs to $ K_n^+$.

Following the analysis performed in dimension 1, we can prove the
backward Kolmogorov formula:
\begin{theorem}
\label{th:d2_feynman} Under the above notations, the numerical
solution $u^n_K$ at time $n$ and in the cell $K$ has the form:
\begin{equation*}
u^n_K = {\mathbb E}_K \bigl[u^0_{K_n} \bigr].
\end{equation*}
\end{theorem}
Due to the strong similarity with the one-dimensional frame, we do
not repeat here the proof.

In what follows, for $K \in {\Tr}$, we denote by ${\mathbb
E}_K^n$ the conditional expectation ${\mathbb E}_K
[\cdot|K_0,\dots,K_n]$.

\subsection{Random Characteristics}

To follow the one-dimensional strategy, we have to associate a
sequence $(X_n)_{n \geq 0}$ of (random) points with each path of
the Markov chain.
In dimension one, the point $X_n$ is defined as the entering
point in the cell $K_n$. Two points may play this role in the higher
dimensional setting: either the barycenter of the entering face
$K_{n-1} \cap K_n$ or the barycenter, with suitable weights, of all
the possible entering faces in the cell $K_n$. For a given cell $K$,
we thus define $x_{K,L}$ as the barycenter of the face $K \cap L$ if
$K$ and $L$ are adjacent: 
\begin{equation}
\label{eq:x}
x_{K,L} = |K \cap L|^{-1} \int_{K \cap L} x dx,
\end{equation}
and $e_K$ as the following convex combination of $(x_{K,J})_{J \in
\Kp}$:
\begin{equation}
\label{eq:e}
e_K = \left(\sum_{J \in \Kp} q_{K,J} \right)^{-1} \sum_{J \in \Kp}
q_{K,J} x_{K,J},
\end{equation}
with
\begin{equation}
\label{eq:q}
\begin{split}
&q_{K,J} = \frac{\langle
a_{K,J},n_{K,J}\rangle \Delta t |K \cap J|}{|K|} \ {\rm for} \ J \in
\Kp, \\
&q_{K,K} = 1 - \sum_{J \in \Kp}\frac{\langle
a_{K,J},n_{K,J}\rangle \Delta t |K \cap J|}{|K|}, \\
&q_{K,J} = 0 \ {\rm for} \ J \in \Tr \setminus \left(K^+ \cup K
\right). 
\end{split}
\end{equation}
We emphasize that the $(q_{K,J})_{J \in {\Tr}}$ are, at least
in a formal way, the weights associated with the scheme for the
velocity $-a$. We will specify this correspondence below.
We also notice that $e_K$ might be outside $K$ if $K$ is not
convex. This does not matter for the analysis.

With these notations at hand, we define the random characteristic
$X_n$ as
\begin{equation}
\label{eq:X}
\begin{split}
&X_0 = e_{K_0}, \\
&X_n = e_{K_{n}} \ {\rm if} \ K_n = K_{n-1}, n \geq 1 \\
&X_n = x_{K_{n-1},K_n} \ {\rm if} \ K_n \not = K_{n-1}, n \geq 1.
\end{split}
\end{equation}

This definition is quite natural. When $n=0$ or $K_{n} = K_{n-1}$, $n
\geq 1$, $X_n$ is chosen as a remarkable point of the cell $K_n$,
independently of the past before $n$. When $K_n \not = K_{n-1}$, the
choice of $X_n$ expresses the jump from $K_{n-1}$ to $K_n$.

\subsection{Green's Formula}

The result given below explains why the barycenters of the faces are
involved in our analysis. The main argument of the proof relies on the
Green formula, which has a crucial role in the whole story, as easily
guessed from the specific form of the transition probabilities.
(See also \cite[Proposition 3.1]{BouGhiPas05}.)
\begin{prop}
\label{prop:greenstokes} Consider a cell $K$. For any point $x_0$ in
the convex envelope of the cell (we say convex envelope because of
$e_K$, defined as a barycenter),
\begin{equation*}
a_K \Delta t = - \sum_{L \in \Km} p_{K,L} \bigl( x_{K,L} - x_0 \bigr)
+ \sum_{J \in \Kp} q_{K,J}
\bigl( x_{K,J} - x_0 \bigr) + O(h \Delta t).
\end{equation*}
\end{prop}
\noindent {\bf Proof.}
For an index $1 \leq i \leq d$, the Green formula, see \cite[Chapter
3, (3.54)]{mizohata}, yields ($a_i$ and $x_i$ stand for the
$i^{\scriptsize{th}}$ coordinate of $a$ and $x$).
\begin{equation}
\label{eq:greenstokes}
\int_K a_i(x) dx = - \int_K (x_i - (x_0)_i) \div(a)(x) dx
+ \sum_{L \sim K} \int_{K \cap L} (x_i-(x_0)_i)
\langle a(x),n_{K,L} \rangle dx.
\end{equation}
The left-hand side is equal to $|K| (a_K)_i$. By the regularity of
$a$, the first term in the right-hand side is bounded by $O(h |K|)$.
Similarly, the last term in the right-hand side writes
\begin{equation*}
\begin{split}
&\sum_{L \sim K} \int_{K \cap L} (x_i-(x_0)_i) \langle a(x),n_{K,L}
\rangle dx \\
&= \sum_{L \sim K} \langle a_{K,L}, n_{K,L} \rangle \int_{K \cap L}
(x_i - (x_0)_i) dx + \sum_{L \sim K} \int_{K \cap L} (x_i-(x_0)_i)
\langle a(x)-a_{K,L},n_{K,L} \rangle dx \\
&= \sum_{L \sim K} \langle a_{K,L},n_{K,L} \rangle \int_{K \cap L}
(x_i - (x_0)_i) dx + O(h^2) \sum_{L \sim K} |K \cap L|.
\end{split}
\end{equation*}
With \eqref{ass:iso} and \eqref{eq:x} at hand, we deduce that
\begin{equation}
\label{eq:green3}
\begin{split}
&\sum_{L \sim K} \int_{K \cap L} (x_i-(x_0)_i) \langle a(x),n_{K,L}
\rangle dx 
\\
&= \sum_{L \sim K} \langle a_{K,L},n_{K,L} \rangle |K
\cap L| \bigl[ x_{K,L} - x_0 \bigr]_i + O(h |K|).
\end{split}
\end{equation}
From \eqref{eq:greenstokes} and \eqref{eq:green3}, we claim
\begin{equation*}
 a_K \Delta t = \sum_{L \sim K} \frac{\langle a_{K,L},n_{K,L} \rangle
\Delta t |K \cap L|}{|K|}
 \bigl( x_{K,L} - x_0 \bigr) + O(h \Delta t).
\end{equation*}
This completes the proof. \qed

The following corollary is the multi-dimensional counterpart of
Proposition \ref{prop:mean_speed}:

\begin{cor}
\label{cor:greenstokes} For any starting cell $K \in {\Tr}$
(so that we work under ${\mathbb P}_K$) and any $n \geq 0$,
\begin{equation*}
{\mathbb E}_K^n \bigl[ X_{n+1} - e_{K_n} \bigr] = - a(X_n) \Delta t +
O(h \Delta t).
\end{equation*}
\end{cor}

As a simple consequence, when the cell $K_n$ has only one possible
entering face (entering means entering for the random characteristic),
as it is the case in dimension 1 with $a > 0$, $e_{K_n} = X_n$ and the
conditional expectation of the mean displacement $X_{n+1} - X_n$
knowing the past before $n$ is driven by $-a$, as stated in
Proposition \ref{prop:mean_speed}.
\vspace{5pt}
\\
\noindent {\bf Proof.} For $K \in {\Tr}$ and $n \geq 0$,
\begin{equation}
\label{cor:gs1}
{\mathbb E}_K^n \bigl[ X_{n+1} - e_{K_n} \bigr] =
\sum_{L \in K_n^-} p_{K_n,L} \bigl( x_{K_n,L} - e_{K_n} \bigr).
\end{equation}
(Indeed, if $K_{n+1}=K_n$, $X_{n+1}=e_{K_n}$.) Applying Proposition
\ref{prop:greenstokes} with $x_0=e_{K_n}$, we obtain
\begin{equation}
\label{cor:gs2}
a_{K_n} \Delta t = - \sum_{L \in K_n^-} p_{K_n,L} \bigl( x_{K_n,L} -
e_{K_n} \bigr) + \sum_{J \in K_n^+} q_{K_n,J} \bigl(x_{K_n,J} - e_{K_n}
\bigr) + O(h \Delta t).
\end{equation}
By the very definition of $e_{K_n}$, the second sum in the above
right-hand side is zero. Identifying \eqref{cor:gs1} and
\eqref{cor:gs2}, we complete the proof. \qed

\subsection{Set of Problems}

Keeping the one-dimensional strategy in mind, we understand that the
whole problem now consists in estimating the gap $e_{K_n} - X_n$ for
$n \geq 1$. (For $n=0$, it is zero.)

As said above, $e_{K_n}-X_n$ vanishes when the cell $K_n$ admits
only one entering face. Unfortunately, there is no hope to obtain a
similar result, or an estimate of the form $e_{K_n}-X_n= O(h \Delta
t)$, for a cell $K_n$ of general shape.

The reason is purely geometric and is well-understood on a
two-dimensional triangular mesh. Except specific cases, the
cardinality of $\Kp$, $K$ being a triangle, is either one or two. (See
Figure \ref{figure:fig1} below.)

\begin{figure}[htb]
\begin{center}
\psfrag{a}{$a$} \psfrag{Xn}{$X_n = x_{K,J_1}$} \psfrag{K}{$K$}
\psfrag{J1}{$J_1$} \psfrag{J2}{$J_2$} \psfrag{J3}{$L$}
\psfrag{e}{$e_K$}
\psfrag{x}{$x_{K,J_2}$}
\psfrag{y}{$x_{K,L}$}
\psfrag{g}{gap}
\includegraphics[
width=0.5\textwidth,angle=0] {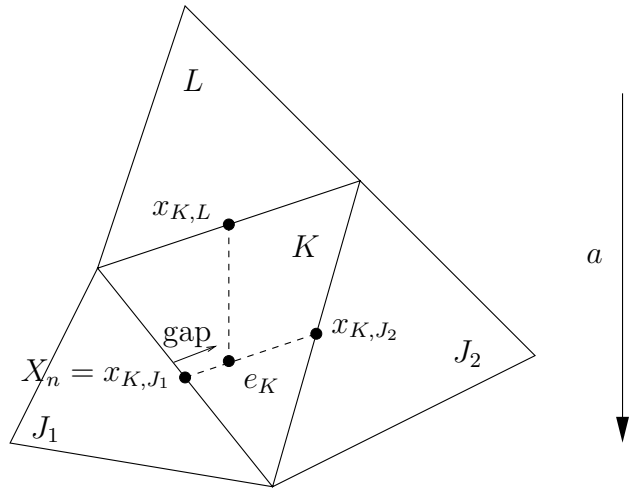} \caption{Example for the
gap $e_{K_n}-X_n$} \label{figure:fig1}
\end{center}
\end{figure}

In this picture, the velocity $a$ is constant, oriented from
the top to the bottom. The triangles $J_1$, $J_2$ and $L$
have only one entering edge and the triangle $K$ has two entering
edges. By Corollary \ref{cor:greenstokes}, the mean displacements in
the triangles $J_1$, $J_2$ and $L$ are driven by $-a$. The story is
different for $K$.

The position of $e_K$ is determined in the following way: $e_K$
belongs to the segment $[x_{K,J_1},x_{K,J_2}]$ and, by Proposition
\ref{prop:greenstokes}, $x_{K,L}-e_K$ is parallel to $a$ (the term
$O(h \Delta t)$ in Proposition \ref{prop:greenstokes} is zero since
$a$ is constant). Assuming that at time $n-1$, $K_{n-1}=J_1$ and
that, at time $n$, $K_n=K$, $X_n$ is the middle of the edge $K \cap
J_1$. We understand that, in this case, $e_{K_n} - X_n$ is of the
same order as $h$.

The key point in our analysis follows from a simple observation:
when reversing the velocity field $a$ in Figure \ref{figure:fig1},
i.e. when changing $a$ into $-a$, the triangles with one entering
edge turn into triangles with two entering edges and the triangle
with two entering edges turns into a triangle with one entering
edge. Said in a very naive way, the bad triangle $K$ for $a$ becomes
a good triangle for $-a$.

Here is another simple fact: the deterministic characteristics for
$-a$ correspond to the characteristics for $a$, but reversed in
time.

These two observations lead to the following idea: in what follows,
we hope to get rid of the remaining gaps $(e_{K_n} - X_n)_{n \geq
1}$ by reversing the random characteristics.

\subsection{Reversing the Markov Chain in the Divergence-Free
Setting.}

As announced above, the weights $(q_{K,J})_{K,J \in {\Tr}}$
in \eqref{eq:q} correspond, at least in a formal sense, to the
scheme associated with the velocity field $-a$. We say ``formally''
because the CFL condition for the field $-a$ may fail: as shown
below, the term $q_{K,K}$ may be negative.

However, when $a$ is divergence-free, the CFL condition holds for $-a$
and $q$ is the analogue of $p$ for the field $-a$:
\begin{prop}
\label{prop:-adiv0} Assume that $a$ is divergence-free. Then,
\begin{equation*}
\forall K \in {\Tr}, \ \sum_{J \sim K} q_{K,J} = \sum_{L \sim
K} p_{K,L},
\end{equation*}
so that $q_{K,K}=p_{K,K}$ for every cell $K$. In particular, the
CFL condition holds for $-a$ and $q$ is a Markovian kernel, i.e.
$q_{K,J} \geq 0$ for all $K,J \in {\Tr}$ and $\sum_{J \in
{\Tr}} q_{K,J} =1$.
\end{prop}
\noindent{\bf Proof.} The proof follows again from Green's
formula. Indeed, for every cell $K$,
\begin{equation}
\label{eq:greenstokesinvmeas}
\int_K \div(a)(x) dx = 0 = \sum_{L \sim K} |K \cap L| \langle
a_{K,L},n_{K,L} \rangle. 
\end{equation}
The terms in the right-hand side are equal to $- p_{K,L} |K|/\Delta
t$ if $L \in \Km$ and $q_{K,L} |K| / \Delta t$ if $L \in
\Kp$. This completes the proof. \qed

In the rest of this section, we assume $a$ to be divergence-free.
The point is to understand the connection between the original
Markov chain $(K_n)_{n \geq 0}$, associated with $a$, and the
Markovian kernel $(q_{K,J})_{K,J \in {\Tr}}$. In what
follows, we show that, under suitable conditions, the random process
obtained by reversing the chain $(K_n)_{n \geq 0}$ is a Markov chain
with $(q_{K,J})_{K,J \in {\Tr}}$ as probability transitions.

Reversing a Markov chain is a standard procedure
in probability theory. (See e.g. \cite[Section 1.9]{norris}.) The
reversed process is always a Markov chain, but the transition
probabilities may be highly non trivial and do depend, in almost every
case, on time.
Anyhow, the law of the reversed chain is easily computable when the
chain is initialized with an invariant probability.

Assume for the moment that the chain $(K_n)_{n \geq 0}$ admits an
invariant probability $\nu$, i.e. there exists a probability $\nu$
on the set of cells such that
\begin{equation*}
\forall K \in {\Tr}, \ \nu(K) = \sum_{J \in {\Tr}}
\nu(J) p_{J,K},
\end{equation*}
and pick the initial cell of the Markov chain randomly, with
respect to the probability measure $\nu$. We then denote the
corresponding probability measure on $(\Omega,{\mathcal A})$ by
${\mathbb P}_{\nu}$ as we did for a deterministic starting cell. Under
${\mathbb P}_{\nu}$, $K_0$ is random and its distribution is given
by $\nu$ itself. The measure ${\mathbb P}_{\nu}$ may be decomposed
along the measures $({\mathbb P}_K)_{K \in {\Tr}}$ according
to the formula ${\mathbb P}_{\nu} = \sum_{K \in {\Tr}} \nu(K)
{\mathbb P}_K$. By invariance of $\nu$, $K_n$ follows, for each $n
\geq 0$, the law $\nu$, i.e.
\begin{equation*}
\forall n \geq 0, \ \forall K \in {\Tr}, \ {\mathbb
P}_{\nu}\{K_n = K\} = \nu(K).
\end{equation*}
(See \cite[Section 1.7]{norris}.) Denoting by
$(\overset{\leftarrow}{K}^{\underset{N}{}}_0,\dots,\lKN_N) =
(K_N,\dots,K_0)$ the reversed chain from time $N$ to time 0, we have,
for $K,J \in {\Tr}$ and $0 \leq n \leq N-1$,
\begin{equation}
\label{eq:invmeas00}
{\mathbb P}_{\nu}\{\lKN_{n+1} = J |\lKN_n = K\}
= \frac{{\mathbb P}_{\nu}\{K_{N-n-1} = J,K_{N-n} = K\}}
{{\mathbb P}_{\nu}\{K_{N-n} = K\}}
= \frac{\nu(J)}{\nu(K)}p_{J,K},
\end{equation}
so that the reversed chain is homogeneous (i.e. the transition
probabilities don't depend on time).

Since the velocity field $a$ is divergence-free, the transport
equation conserves the mass with respect to the Lebesgue measure. In
this framework, the Lebesgue measure is invariant for the Markov chain
$(K_n)_{n \geq 0}$:
\begin{prop}
\label{prop:invmeas} Assume that $a$ is divergence-free. Then,
the Lebesgue measure is invariant for the
Markov chain, i.e.
\begin{equation*}
\forall K \in {\Tr}, \ |K| = \sum_{J \in {\Tr}} |J|
p_{J,K}.
\end{equation*}
\end{prop}
\noindent {\bf Proof.} For a given cell $K$, we have
\begin{equation}
\label{eq:invmeas1}
\begin{split}
\sum_{J \in {\Tr}} |J| p_{J,K} &= \sum_{J \in \Kp}
\langle a_{J,K},n_{J,K} \rangle \Delta t |J \cap K|
 + |K| + \sum_{J \in \Km}\langle a_{J,K},n_{J,K}
 \rangle \Delta t |J \cap K|.
\end{split}
\end{equation}
By \eqref{eq:greenstokesinvmeas},
\begin{equation*}
\sum_{J \sim K} |K \cap J| \langle a_{J,K},n_{J,K} \rangle =0,
\end{equation*}
since $a$ is divergence-free. Returning to \eqref{eq:invmeas1}, we
obtain $ \sum_{J \in {\Tr}} |J| p_{J,K} = |K|$. \qed

\subsection{Analysis in the Divergence-Free and Periodic Setting}

The Lebesgue measure is not of finite mass on $\R^d$, but it is a
probability measure on the torus $\R^d/{\mathbb Z}^d$.

To explain how we make use of time reversal, we thus assume, for the
moment, that the problem is periodic, of period one in each
direction. (Since the transport is of finite speed, this is not a big
deal. Anyhow, the proofs of the main results of the paper, given in
the next section, are performed without any periodicity
assumption. The current paragraph is purely pedagogical.)
This means that both the velocity $a$ and the mesh are periodic, of
period one in each direction of the space.

As a consequence, we can see the Markov chain $(K_n)_{n \geq 0}$ as a
Markov chain with values in the set ${\Tr}/{\mathbb Z}^d$,
i.e. in the space of classes of cells for the equivalence relation
induced by the periodicity. The probability of jumping from the
class of the cell $K$ to the class of the cell $L$ is given by the
rate $p_{K,L}$: by periodicity, this rate doesn't depend on the
choices of $K$ and $L$.

\begin{prop}
\label{prop:invmeastorus} Assume that $a$ is divergence-free and
that both $a$ and the mesh are periodic of period one in each
direction of the space, then the Lebesgue measure on the torus
induces an invariant probability for the Markov chain, i.e.
\begin{equation*}
\forall K \in {\Tr}/{\mathbb Z}^d, \ |K| = \sum_{J \in
{\Tr}/{\mathbb Z}^d, \ J \sim K} |J| p_{J,K}.
\end{equation*}
(For $J \in {\Tr}/{\mathbb Z}^d$, we denote by $|J|$ the
common volume of all the cells of $J$.)

In particular, by \eqref{eq:invmeas00}, the reversed chain has the
following transition probabilities
\begin{equation*}
{\mathbb P}_{\mu} \bigl\{ \lKN_{n+1} = J |\lKN_n =K \bigr\} =
q_{K,J}.
\end{equation*}
(In the above equality, $\mu$ denotes the Lebesgue measure on the
torus. Under ${\mathbb P}_{\mu}$, the starting class of cells of the
chain is chosen according to the Lebesgue measure.)
\end{prop}
\noindent {\bf Proof.} To check the last assertion, we note that,
for two adjacent cells $J$ and $K$, $p_{J,K} >0$ if and only if $K
\in J^-$, that is $J \in \Kp$. Hence, $p_{J,K}
>0 \Leftrightarrow q_{K,J} >0$ and in this case
\begin{equation}
 \label{eq:reverse100}
q_{K,J} = \frac{|J|}{|K|} p_{J,K},
\end{equation}
so that the relation is always true (even if one term vanishes).
\qed

We recall the interpretation of Proposition \ref{prop:invmeastorus}:
the reversed chain (seen as a chain with values in $\Tr /{\mathbb
Z}^d$), when the chain $(K_n)_{n \geq 0}$ is initialized with the
Lebesgue measure, is nothing but the chain associated with $-a$.

Before applying Proposition \ref{prop:invmeastorus} to the analysis of
the numerical scheme, we have to specify the construction of the random
characteristic $(X_n)_{n \geq 0}$ in the periodic setting. In this
case, $(X_n)_{n \geq 0}$ is seen as a path in $\R^d$ and not in
$\R^d/{\mathbb Z}^d$. This amounts to consider a sequence of
representatives $(\hat{K}_n)_{n \geq 0}$ for the chain $(K_n)_{n
\geq 0}$. We choose $\hat{K}_0$ as the only representative $K^0_0$
of $K_0$ such that $e_{K_0^0}$ has coordinates in $[0,1)$, that is
$\hat{K}_0=K_0^0$ (and $X_0 = e_{K_0^0}$). For any $n \geq 0$,
$\hat{K}_{n+1}$ is the unique representative of $K_{n+1}$ such that
$\hat{K}_{n+1} \in (\hat{K}_n)^-$. With this sequence of
representatives at hand, we can build up the sequence $(X_n)_{n \geq
0}$ according to \eqref{eq:X}. We are now in position to estimate the
gap $e_{\hat{K}_n} - X_n$, $n \geq 1$:

\begin{prop}\label{prop:etorus}
Under the assumptions of Proposition \ref{prop:invmeastorus}, for any
$N \geq 0$ and any $n \in \{0,\dots,N\}$,
\begin{equation*}
{\mathbb E}_{\mu}\bigl[ e_{\hat{K}_n} - X_n |K_n,\dots,K_N \bigr] = 0,
\end{equation*}
${\mathbb E}_{\mu}$ denoting the expectation under ${\mathbb
P}_{\mu}$. In particular, under ${\mathbb P}_{\mu}$,
the sequence $(\lMN_n)_{0 \leq n \leq N}$, given by
\begin{equation*}
\lMN_0 = \lMN_1=0, \ \lMN_n = \sum_{k=N-n+1}^{N-1} \bigl(
e_{\hat{K}_{k}} - X_{k} \bigr), \ 2 \leq n \leq N,
\end{equation*}
is a martingale for the backward
filtration $(\sigma(K_{N-n},\dots,K_N)=\sigma(\lKN_0,\dots,\lKN_n))_{0
\leq n \leq N}$.
\end{prop}
\noindent
{\bf Proof.} We can assume that $1 \leq n \leq N$, since
$e_{\hat{K}_0} - X_0=0$. We then emphasize that $\hat{K}_n$ isn't
measurable with respect to $\sigma(K_n,\dots,K_N)$. This is the main
difficulty of the proof. Indeed, $\hat{K}_n$ depends both on $K_n$ and
on the initial representative of $K_0$. Anyhow, the difference
$e_{\hat{K}_n} - X_n$ doesn't depend on the representatives chosen for
$K_{n-1}$ and $K_n$.

Indeed, recalling that $K_n^0$ is the unique representative of
$K_n$ such that $e_{K_n^0}$ has coordinates in $[0,1)$ , we can
always write
\begin{equation}
\label{eq:ehat} e_{\hat{K}_n} - X_n = \sum_{J \in (\hat{K}_n)^+}
\bigl( e_{\hat{K}_n} - x_{\hat{K}_n,J} \bigr) {\mathbf
1}_{\{\hat{K}_{n-1} = J\}} = \sum_{J \in (K_n^0)^+} \bigl( e_{K_n^0}
- x_{K_n^0,J} \bigr) {\mathbf 1}_{\{K_{n-1} \ni J\}}.
\end{equation}
(Above, $K_{n-1} \ni J$ means that $J$ is a representative of
$K_{n-1}$.) Hence,
\begin{equation*}
{\mathbb E}_{\mu} \bigl[ e_{\hat{K}_n} - X_n | K_n,\dots,K_N \bigr]
= \sum_{J \in (K_n^0)^+} \bigl( e_{K_n^0} - x_{K_n^0,J} \bigr)
{\mathbb P}_{\mu} \bigl\{ K_{n-1} \ni J | K_n,\dots,K_N \bigr\}.
\end{equation*}
The probability ${\mathbb P}_{\mu} \{ K_{n-1} \ni J | K_n,\dots,K_N
\}$ is also ${\mathbb P}_{\mu} \{ \lKN_{N-n+1} \ni J |
\lKN_0,\dots,\lKN_{N-n} \}$. By Proposition \ref{prop:invmeastorus},
it is equal to $q_{K_n^0,J}$. By \eqref{eq:e}, we deduce that
${\mathbb E}_{\mu} [ e_{\hat{K}_n} - X_n | K_n,\dots,K_N ]=0$. By
\eqref{eq:ehat}, we know that $e_{\hat{K}_{n+1}} - X_{n+1}$ is
$\sigma(K_n,\dots,K_N)$-measurable for $0 \leq n \leq N-1$. We
easily deduce the martingale property. \qed

Following \eqref{d1:300}, \eqref{eq:4}, \eqref{eq:varianced1}, we
complete the analysis by estimating, under ${\mathbb P}_{\mu}$, the
fluctuations of the random characteristic around the deterministic
characteristic. 
\begin{prop}
\label{prop:estimtorus} There exists a constant $C \geq 0$ such
that, for any $N \geq 1$,
\begin{equation*}
{\mathbb E}_{\mu} \bigl[ |X_N - Z(N \Delta t,X_0)| \bigr] \leq
C \bigl[ (N h \Delta t)^{1/2} + N h \Delta t + N \Delta t^2 \bigr]
\exp(\kappa N \Delta t).
\end{equation*}
\end{prop}
\noindent {\bf Proof.} Following the proof of \eqref{d1:300} and
keeping the identity $e_{\hat{K}_0} - X_{0}=0$  in mind, it is
sufficient to focus on
\begin{equation*}
X_N - X_0 + \sum_{k=0}^{N-1} a(X_{k}) \Delta t = \sum_{k=0}^{N-1}
(X_{k+1}-e_{\hat{K}_k} + a(X_{k}) \Delta t) + \sum_{k=1}^{N-1}
(e_{\hat{K}_k} - X_{k}).
\end{equation*}
It is plain to see that Corollary \ref{cor:greenstokes} is still true
under ${\mathbb P}_{\mu}$ (and not under ${\mathbb P}_K$), so
that
\begin{equation*} 
X_N - X_0 + \sum_{k=0}^{N-1} a(X_{k}) \Delta t =
\sum_{k=0}^{N-1} \bigl(X_{k+1}-e_{\hat{K}_k} - {\mathbb
E}_{\mu}^k(X_{k+1}-e_{\hat{K}_k}) \bigr) + \sum_{k=1}^{N-1}
(e_{\hat{K}_k} - X_{k}) + O(N h \Delta t).
\end{equation*}
(Above, ${\mathbb E}_{\mu}^k[\cdot] = {\mathbb E}_{\mu}[\cdot|K_0,
\dots,K_k]$.) Setting, for all $0 \leq n \leq N$, $M_n =
\sum_{k=0}^{n-1} (X_{k+1}-e_{\hat{K}_k} - {\mathbb
E}_{\mu}^k(X_{k+1}-X_k))$ (with $M_0=0$), and using the notation
introduced in Proposition \ref{prop:etorus}, we write
\begin{equation*}
X_N - X_0 + \sum_{k=0}^{N-1} a(X_{k}) \Delta t = M_N + \lMN_N + O(N h
\Delta t).
\end{equation*}
We let the reader check that $(M_n)_{0 \leq n \leq N}$ is a martingale
with respect to the (forward) filtration $(\sigma(K_0,\dots,K_n))_{0
\leq n \leq N}$. By Proposition \ref{prop:etorus}, $(\lMN_n)_{0 \leq n
\leq N}$ is a martingale with respect to the backward filtration.

Following \eqref{eq:4}, we have
\begin{equation*}
{\mathbb E}_{\mu} \bigl[|M_N|^2 \bigr] = \sum_{n=0}^{N-1} {\mathbb
E}_{\mu} \bigl[ |M_{n+1} - M_n |^2 \bigr], \ {\mathbb E}_{\mu}
\bigl[\bigl| \lMN_N \bigr|^2 \bigr] = \sum_{n=0}^{N-1} {\mathbb
E}_{\mu} \bigl[ |\lMN_{n+1} - \lMN_n |^2 \bigr].
\end{equation*}

Following \eqref{eq:varianced1}, we have, for all $0 \leq n \leq N-1$,
\begin{equation}
\label{eq:condivar}
\begin{split}
{\mathbb E}_{\mu} \bigl[ |M_{n+1} - M_n|^2 \bigr]
&\leq {\mathbb E}_{\mu} \bigl[ |X_{n+1} - e_{\hat{K}_n}|^2 \bigr] \\
&\leq {\mathbb E}_{\mu} \bigl[ \sum_{L \in (\hat{K}_n)^-}
p_{\hat{K}_n,L} |x_{\hat{K}_n,L} - e_{\hat{K}_n}|^2 \bigr] \\
&\leq h^2 \sup_{K} \sum_{L \sim K} p_{K,L} \leq \|a\|_{\infty}
\sup_{K} \bigl[ |K|^{-1} \sum_{L \sim K} |K \cap L| \bigr] h^2 \Delta
t.
\end{split}
\end{equation}
By \eqref{ass:iso}, we deduce ${\mathbb E}_{\mu} [|M_{n+1} - M_n
|^2] \leq \alpha \|a\|_{\infty} h \Delta t$. By a similar argument,
we obtain ${\mathbb E}_{\mu} [|\lMN_{n+1} - \lMN_n |^2] \leq \alpha
\|a\|_{\infty} h \Delta t$. We deduce
\begin{equation*}
{\mathbb E}_{\mu} \bigl[|M_N|^2 \bigr] = O(hN \Delta t), \ {\mathbb
E}_{\mu} \bigl[\bigl| \lMN_N \bigr|^2 \bigr] = O(h N \Delta t).
\end{equation*}
We then complete the proof as in the one-dimensional case.
\qed

\subsection{$L^1$-Error in the Divergence-Free and Periodic Setting
with $u^0$ Lipschitz continuous}

As a by-product, we obtain the following estimate for the error of the
numerical scheme when $u^0$ is Lipschitz continuous:

\begin{theorem}\label{cor:L1conv}
Assume that the hypotheses of the beginning of Section \ref{sec:2} are
satisfied. Assume moreover that $a$ is divergence-free and that both
the velocity $a$ and the mesh are periodic of period one in each
direction of the space. Assume also that $u^0$ is $\kappa$-Lipschitz
continuous. Then, there exists a constant $C \geq 0$ such that
\begin{equation*}
\sum_{K\in{\Tr}/{\mathbb Z}^d} \int_{K} |u^{N}_K - u(N \Delta
t,x)| dx
 \leq C \bigl(
 (h N \Delta t)^{1/2} + h N \Delta t + N \Delta t^2\bigr)
 \exp(\kappa N \Delta t).
\end{equation*}
\end{theorem}
\noindent {\bf Proof.} Consider a cell $K$. We know that $u^N_K$ is
nothing but
\begin{equation*}
u_K^N = {\mathbb E}_K \bigl[ u^0_{K_N} \bigr] = {\mathbb E}_K \bigl[
u^0(X_N)\bigr] +O(h).
\end{equation*}
Moreover, for all $x \in
K$, $|X_0-x| \leq h$ with probability one under ${\mathbb P}_{K}$.
By stability of the solutions of \eqref{carac}, $|Z(N \Delta
t,x) - Z(N \Delta t,X_0)| \leq h \exp(\kappa N \Delta t)$ under
${\mathbb P}_{K}$. We deduce
\begin{equation*}
\forall x \in K, \ u(N\Delta t,x) = u^0\bigl(Z(N \Delta t,x)\bigr) =
{\mathbb E}_K \bigl[ u^0 \bigl(Z(N \Delta t,X_0) \bigr) \bigr] +
O\bigl(h \exp(\kappa N \Delta t) \bigr).
\end{equation*}
Hence,
\begin{equation*}
\begin{split}
\sum_{K \in{\Tr}/{\mathbb Z}^d} \int_{K} |u^{N}_K - u(N
\Delta t,x)| dx &\leq \kappa \sum_{K\in{\Tr}/{\mathbb Z}^d}
|K| {\mathbb E}_K \bigl[ |X_N - Z(N \Delta t,X_0)| \bigr] +
O\bigl(h \exp(\kappa N \Delta t) \bigr) \\
&\leq \kappa {\mathbb E}_{\mu} \bigl[ |X_N - Z(N \Delta t,X_0)|
\bigr] + O\bigl(h \exp(\kappa N \Delta t) \bigr).
\end{split}
\end{equation*}
This completes the proof.\qed

\section{Analysis in the General Setting}\label{gensec}

We now turn to the general case and analyze the error of the
numerical scheme, both in the $L^1$ sense and in the $L^{\infty}$
sense, $u^0$ being respectively of bounded variation and Lipschitz
continuous. We thus forget the periodic setting and the
divergence-free condition. Anyhow, for technical reasons explained
below, the time step $\Delta t$ is required to be small when the
divergence of $a$ is large, i.e.
\begin{equation}\label{highd:deltat}
\exists \eta \in (0,1), \ \|\div(a)\|_{\infty} \Delta t < 1-\eta.
\end{equation}

We keep the notations of Section \ref{sec:2}. To simplify the form
of the final bounds, we assume that $h \leq 1$ and that $\Delta t
\leq \theta h$ for some $\theta >0$. The constants $C$, $C'$
and $c$ below may depend on the parameters specified in Section
\ref{sec:2}, on $\eta$ and on $\theta$. The values of these
``constants'' may vary from line to line.

\subsection{Strategy}

In dimension one, we were able to analyze the error in the
$L^{\infty}$ sense by investigating the distance between $X_N$ and
$X_0 - \Delta t \sum_{n=0}^{N-1} a(X_n)$ for any arbitrary initial
cell $K$ of the mesh. In the previous section, the result was given
in the $L^1$ norm since the chain was initialized with the Lebesgue
measure. In what follows, the approach is halfway.

The idea is the following. We pick up the starting cell of the chain
with respect to the Lebesgue measure among the cells included in a
ball of radius $h^{1/2}$ and centered at the origin (or at any other
arbitrarily prescribed point). The bounds we then obtain for the
error of the numerical scheme hold in the $L^1$ sense, but locally in
the ball. In other words, we manage to bound
\begin{equation*}
h^{-d/2} \sup_{x \in \R^d} \int_{B(x,h^{1/2})} |u^N(y) - u(N \Delta
t,y)| dy
\end{equation*}
by a constant times $h^{1/2}$ (up to remaining terms in $N$ and
$\Delta t$), $u^0$ being Lipschitz continuous. (In the above
expression, $u^N(y)$ stands for $u^N_K$ when $y$ belongs to the
interior of $K$.) By an approximation procedure, we deduce that the
scheme is of order $1/2$ for the (global) $L^1$ norm when $u^0$ is
of bounded variation.

Actually, we can perform the same analysis by replacing the local
$L^1$ norm by the local $L^p$ norm, $p$ being greater than one. We
then derive that, for any small positive $\varepsilon$, the scheme
is of order $1/2-\varepsilon$ for the $L^{\infty}$ norm when $u^0$
is Lipschitz continuous. Moreover, by translation, it is sufficient to
prove this estimate for $x = 0$. 

For all these reasons, the quantities of interest are
\begin{equation*}
{\mathcal Q}_p^N = h^{d/2} \sum_{K \in {\Tr}_0} {\mathbb
E}_{K} \bigl[ \bigl|X_N-X_0 + \Delta t \sum_{n=0}^{N-1} a(X_n)
\bigr|^p \bigr], \quad N, p \geq 1,
\end{equation*}
where ${\Tr}_0$ stands for the set of cells $K$ such that
$|e_K| \leq h^{1/2}$. By \eqref{ass:iso}, the cardinality of
${\Tr}_0$ is bounded by $Ch^{-d/2}$ for some positive
constant $C$.

For given $N,p \geq 1$, we are going to decompose ${\mathcal Q}_p^N$
along all the possible paths of the chain. To do so, we distinguish
the paths according to their fluctuations around the velocity field
$-a$. We write
\begin{equation*}
{\mathcal Q}_p^N \leq h^{(d+p)/2} \sum_{K \in {\Tr}_0}
\sum_{k \geq 0} (k+1)^p {\mathbb P}_K \bigl\{k h^{1/2} \leq \sup_{1
\leq n \leq N} \bigl| X_n - X_0 + \Delta t \sum_{i=0}^{n-1} a(X_i)
\bigr| < (k+1)h^{1/2} \bigr\}.
\end{equation*}

To estimate the above right-hand side, we will use a time reversal
argument, as in Section \ref{sec:2}. Therefore, we are mainly
interested in the location of the arrival cell $K_N$ on the event
$\{k h^{1/2} \leq \sup_{1 \leq n \leq N} | X_n - X_0 + \Delta t
\sum_{i=0}^{n-1} a(X_i) | < (k+1)h^{1/2}\}$. Thus, for any $k \geq
0$, we denote by ${\Tr}_k^N$ the set of cells $J_N$ such
that there exists an $N$-tuple $(J_0,\dots,J_{N-1})$, with $J_0 \in
{\Tr}_0$, satisfying, for all $i \in \{0,\dots,N-1\}$, either
$J_{i+1} = J_i$ or
 $J_{i+1} \in J_{i}^-$, and
\begin{equation}\label{eq:general10000}
k h^{1/2} \leq \sup_{1 \leq n \leq N} |y_{J_{n-1},J_n} - e_{J_0} +
\Delta t \sum_{i=0}^{n-1} a(y_{J_{i-1},J_{i}}) | < (k+1) h^{1/2},
\end{equation}
with $y_{J_i,J_{i+1}} = e_{J_i}$if $J_i=J_{i+1}$ and
$y_{J_i,J_{i+1}}=x_{J_i,J_{i+1}}$ if $J_{i+1} \in J_i^-$ (and
$y_{J_{-1},J_0}=e_{J_0}$). Under ${\mathbb P}_K$, for $K \in{\Tr}_0$,
we have $K_N\in \Tr_k^N$ on the event $\{k h^{1/2} \leq \sup_{1 \leq n
\leq N} | X_n - X_0 + \Delta t \sum_{i=0}^{n-1} a(X_i) |<
(k+1)h^{1/2}\}$. 

Thus,
\begin{equation}
\label{eq:Q}
\begin{split}
{\mathcal Q}_p^N &\leq h^{(d+p)/2}\sum_{k \geq 0} \sum_{K \in
{\Tr}_0} \sum_{L \in {\Tr}_k^N} \bigl[ (k+1)^p
\\
&\hspace{15pt} \times {\mathbb P}_K \bigl\{k h^{1/2} \leq \sup_{1
\leq n \leq N} \bigl| X_n - X_0 + \Delta t \sum_{i=0}^{n-1} a(X_i)
\bigr| < (k+1)h^{1/2},K_N=L \bigr\} \bigr].
\end{split}
\end{equation}

In the next lemma, we estimate the cardinality of
$\Tr^N_k$. Because of the translation action of the transport
equation, $\sharp[ {\Tr}^N_k]$ is of the same order as $\sharp[
{\Tr}_0]$:

\begin{lem}
\label{lem:general40} There exists a constant $C>0$ such that, for
all $k \geq 0$,
\begin{equation*}
\sharp [ {\Tr}^N_k ] \leq C \exp(C N \Delta t) (k+1)^d
h^{-d/2}.
\end{equation*}
\end{lem}
\noindent {\bf Proof.} We fix $k \geq 0$ and we consider a sequence
of cells $(J_0,\dots,J_N)$, $J_0 \in {\Tr}_0$, satisfying
\eqref{eq:general10000}. Setting $(y^{0},\dots,y^{N}) =
({e}_{J_0},y_{J_0,J_1},\dots,y_{J_{N-1},J_N})$,
\begin{equation*}
\sup_{1 \leq n \leq N} \bigl| y^{n} - y^{0} + \Delta t
\sum_{i=0}^{n-1} a(y^i) \bigr| \leq (k+1) h^{1/2}.
\end{equation*}
Plugging the characteristic of the transport equation (see
\eqref{carac}), we deduce (recall that $h \leq 1$ and $\Delta t
\leq \theta h$ by assumption)
\begin{equation*}
\begin{split}
&\sup_{1 \leq n \leq N} \bigl|y^{n} - Z(n \Delta t,y^{0}) + \Delta
t \sum_{i=0}^{n-1} \bigl[ a(y^i) - a \bigl(Z(i \Delta t,y^{0})
\bigr) \bigr] \bigr|
\\
&\hspace{15pt} \leq (k+1)h^{1/2} + C N \Delta t^2
\leq (k+1)h^{1/2} ( 1 + C N \Delta t),
\end{split}
\end{equation*}
for some constant $C>0$. By the Lipschitz property of $a$ and the
Gronwall lemma, it is plain to deduce (up to a new value of $C$)
\begin{equation*}
\bigl|y^{N} - Z(N \Delta t,y^{0}) \bigr| \leq C \exp( C N \Delta
t) (k+1) h^{1/2}.
\end{equation*}
Since $|y^{0} | \leq h^{1/2}$ ($J_0 \in {\Tr}_0$), we deduce,
by stability of the solutions to \eqref{carac}, that
\begin{equation*}
\bigl|y^{N} - Z(N \Delta t,0) \bigr| \leq C \exp( C N \Delta t)
(k+1) h^{1/2}.
\end{equation*}
Every point $x \in J_N$ satisfies the same property since the
diameter of $J_N$ is bounded by $h$. We deduce that there exists a
constant $C$ such that $J_N$ is included in the ball of center $Z(N
\Delta t,0)$ and of radius $C \exp (C N \Delta t) (k+1) h^{1/2}$. In
other words, all the cells in ${\Tr}^N_k$ are included in
this ball. Up to a modification of $C$, the volume of the ball is $C
\exp(C N \Delta t) (k+1)^d h^{d/2}$. Since the volume of a given
cell is greater than $\beta^{-1} h^d$ (see Assumption
\eqref{ass:iso}), the cardinality of ${\Tr}_k^N$ is
bounded by $C \exp(C N \Delta t) (k+1)^d h^{-d/2}$ for a new value
of the constant $C$. \qed

\subsection{Application of Section \ref{sec:2}}

In light of Section \ref{sec:2}, we introduce the following
decomposition
\begin{equation*}
X_N-X_0 + \Delta t \sum_{n=0}^{N-1} a(X_n) = S_N + R_N,
\end{equation*}
with
\begin{equation*}
S_0=R_0=0, \ S_n = \sum_{i=0}^{n-1} \bigl[X_{i+1} - e_{K_i} + \Delta
t \, a(X_i) \bigr], \ R_n = \sum_{i=0}^{n-1} \bigl[ e_{K_i} - X_i
\bigr], \ n\geq 1.
\end{equation*}
On the event $\{k h^{1/2} \leq \sup_{1 \leq n \leq N} | X_n - X_0 +
\Delta t \sum_{i=0}^{n-1} a(X_i) | < (k+1)h^{1/2}\}$, we have
$\sup_{1 \leq n \leq N} |S_n| \geq (k/2)h^{1/2}$ or $\sup_{1 \leq n
\leq N} |R_n| \geq (k/2)h^{1/2}$. By \eqref{eq:Q}, we obtain
\begin{equation}
\label{eq:R+S}
\begin{split}
h^{-(d+p)/2} {\mathcal Q}_p^N &\leq \sum_{k \geq 0} \sum_{K \in
{\Tr}_0} \sum_{L \in {\Tr}_k^N} (k+1)^p {\mathbb
P}_K \bigl\{ \sup_{0 \leq n \leq N} |S_n| \geq \frac{k}{2} h^{1/2},
K_N=L \bigr\}
\\
&\hspace{15pt} + \sum_{k \geq 0} \sum_{K \in {\Tr}_0}
\sum_{L \in {\Tr}_k^N} (k+1)^p {\mathbb P}_K \bigl\{
\sup_{0 \leq n \leq N} |R_n| \geq \frac{k}{2} h^{1/2}, K_N
=L\bigr\}.
\end{split}
\end{equation}

We wish to apply Corollary \ref{cor:greenstokes} to treat the first
term in the above right-hand side. Since $\sharp[{\Tr}_0]
\leq C h^{-d/2}$, we have
\begin{equation}
\label{eq:Ssup}
\begin{split}
&\sum_{k \geq 0} \sum_{K \in {\Tr}_0} \sum_{L \in \Tr_k^N} (k+1)^p
{\mathbb P}_K \bigl\{ \sup_{0 \leq n \leq N} |S_n| \geq \frac{k}{2}
h^{1/2}, K_N=L \bigr\} \\
&\hspace{5pt} \leq C h^{-d/2} \sum_{k \geq 0} (k+1)^p \sup_{K}
{\mathbb P}_K \bigl\{ \sup_{0 \leq n \leq N} |S_n| \geq \frac{k}{2}
h^{1/2} \bigr\}
\\
&\hspace{5pt} \leq C h^{-d/2} \biggl[ h^{1/2} N \Delta t (1+ C N
\Delta t)^p + \sum_{k \geq 0} (k+1)^p
 \bigl[ \exp \bigl( - \frac{k^2}{C N \Delta t} \bigr)
+ \exp \bigl( - \frac{k}{C h^{1/2}} \bigr) \bigr] \biggr],
\end{split}
\end{equation}
the last line following from
\begin{lem}
\label{lem:S} There exists a constant $C >0$ such that, for $k > C N
\Delta t \, h^{1/2}$,
\begin{equation*}
 \sup_{K \in {\mathcal T}} {\mathbb P}_K \bigl\{ \sup_{0 \leq
n \leq N} |S_n| \geq \frac{k}{2} h^{1/2} \bigr\} \leq
 C \bigl[ \exp \bigl( - \frac{k^2}{C N \Delta t} \bigr)
+ \exp \bigl( - \frac{k}{C h^{1/2}} \bigr) \bigr].
\end{equation*}
\end{lem}
\noindent {\bf Proof of Lemma \ref{lem:S}.} We fix the starting cell
$K \in {\Tr}$. (We thus work under ${\mathbb P}_K$.)
Following the proof of Proposition \ref{prop:estimtorus}, we
introduce the sequence
\begin{equation*}
M_0=0, \ M_n = \sum_{i=0}^{n-1} \bigl( X_{i+1} - e_{K_i} -{\mathbb
E}_K^i[X_{i+1} - e_{K_i}] \bigr), \ n \geq 1.
\end{equation*}
It is a martingale with respect to the filtration
$(\sigma(K_0,\dots,K_n))_{n \geq 0}$ (under ${\mathbb P}_K$). By
Corollary \ref{cor:greenstokes}, we have
\begin{equation*}
\sup_{0\leq n \leq N} |M_n - S_n| \leq C N \Delta t \, h,
\end{equation*}
for a positive constant $C$. Thus, on the event $\{ \sup_{0 \leq n
\leq N} |S_n| \geq (k/2)h^{1/2} \}$, $\sup_{0\leq n \leq N} |M_n|
\geq (k/4) h^{1/2}$ or $CN h \Delta t \geq (k/4) h^{1/2}$. The
latter is impossible if $k > 4CN h^{1/2}\Delta t $. We deduce that
\begin{equation*}
{\mathbb P}_K \bigl\{ \sup_{0 \leq n \leq N} |S_n| \geq \frac{k}{2}
h^{1/2} \bigr\} \leq {\mathbb P}_K \bigl\{ \sup_{0 \leq n \leq N}
|M_n| \geq \frac{k}{4} h^{1/2} \bigr\},
\end{equation*}
for $k > 4 CN \Delta t \, h^{1/2}$. As in \eqref{eq:condivar}, the
conditional variance of the martingale $(M_n)_{n \geq 0}$ may be
bounded by $C \Delta t \, h$ for a possibly new value of $C$:
\begin{equation*}
\forall n \geq 0, \ {\mathbb E}_K^n \bigl[ |M_{n+1} - M_n|^2
\bigr]\leq C \Delta t \, h.
\end{equation*}
Moreover, the jumps of the martingale are bounded by $h$, i.e.
$|M_{n+1} - M_n| \leq h$ for all $n \geq 0$. Applying Proposition
\ref{prop:concentration} given in Annex to $(h^{-1}M_n)_{n \geq
0}$, we obtain
\begin{equation*}
{\mathbb P}_K \bigl\{ \sup_{0 \leq n \leq N} |M_n| \geq \frac{k}{4}
h^{1/2} \bigr\} \leq C \bigl[ \exp \bigl( - \frac{k^2}{C N \Delta t}
\bigr) + \exp \bigl( - \frac{k}{C h^{1/2}} \bigr) \bigr].
\end{equation*}
(Pay attention, the value of $v$ in Proposition
\ref{prop:concentration} is $v=C \Delta t \, h^{-1}$ since we divide
$(M_n)_{n \geq 0}$ by $h$.)

\subsection{Time Reversal}

To treat the gap term $R_N$ in \eqref{eq:R+S}, we are to reverse the
chain $(K_n)_{n \geq 0}$ as done in Section \ref{sec:2} and then to
compare the law of the reversed chain with the law of the chain
associated with $-a$.

Because of \eqref{eq:greenstokesinvmeas}, we emphasize that the
CFL condition may fail for $-a$, so that we cannot associate a
Markov chain with the weights $(q_{K,J})_{K,J \in {\Tr}}$.
The following proposition says how to modify them to obtain a
Markovian kernel:
\begin{prop}
\label{prop:-adivnot0} Set
\begin{equation*}
\forall K \in {\Tr}, \
\delta_K = |K|^{-1} \int_{K} \div(a)(x)dx.
\end{equation*}
Then, under Condition \eqref{highd:deltat}, the kernel
\begin{equation*}
\begin{split}
&\gamma_{K,J} = (1 + \delta_K \Delta t)^{-1} q_{K,J} \ {\rm for} \
J \in \Kp, \\
&\gamma_{K,K} = 1 - \sum_{J \in \Kp} (1 + \delta_K
\Delta t)^{-1} q_{K,J}, \\
&\gamma_{K,J} =0 \ {\rm for} \ J \in \Tr \setminus \left( \Kp \cup K
\right), 
\end{split}
\end{equation*}
satisfies $\gamma_{K,K} = (1 + \delta_K \Delta t)^{-1} p_{K,K}$. In
particular, it is Markovian, i.e. $\gamma_{K,J} \geq 0$ for all $K,J
\in {\Tr}$ and $\sum_{J \in {\Tr}} \gamma_{K,J}=1$.
\end{prop}
\noindent
By \eqref{highd:deltat}, we emphasize that $\gamma_{K,J}\leq
\eta^{-1} q_{K,J}$ for $J \in \Kp$ and $\gamma_{K,K} \leq \eta^{-1}
p_{K,K}$.
\vspace{5pt}

\noindent{\bf Proof.} Consider a cell $K$. By \eqref{highd:deltat}, we
have $1 + \delta_K \Delta t>0$. We compute $\gamma_{K,K}$. By
integrating by parts the expression of $\delta_K$, as in
\eqref{eq:greenstokesinvmeas}, we get 
\begin{equation*}
- \sum_{L \sim K} p_{K,L} + \sum_{J \sim K} q_{K,J} = \delta_K
\Delta t,
\end{equation*}
that is $p_{K,K} -1 + (1+ \delta_K \Delta t)(1-\gamma_{K,K})=
\delta_K \Delta t$.
We deduce that $\gamma_{K,K} = (1+\delta_K \Delta t)^{-1} p_{K,K}$.
\qed
\vspace{5pt}

The chain associated with the kernel $\gamma$ is denoted by the pair
$(({\mathbb Q}_K)_{K \in {\Tr}},(\Gamma_n)_{n \geq 0})$:
$(\Gamma_n)_{n \geq 0}$ is a sequence of measurable mappings from
$(\Omega,{\mathcal A})$ into the set of cells and $({\mathbb
Q}_K)_{K \in {\Tr}}$ is a family of probability measures on
$(\Omega,{\mathcal A})$, such that $(\Gamma_n)_{n \geq 0}$ is, under
${\mathbb Q}_K$, a Markov chain with $K$ as initial condition and
$\gamma$ as kernel. The expectation under ${\mathbb Q}_K$ is denoted
by ${\mathbb E}^{\mathbb Q}_K$.

The link between the chain $(\Gamma_n)_{n \geq 0}$ and the reversed
chain $(\lKN_n = K_{N-n})_{0 \leq n \leq N}$ is given by

\begin{prop}
\label{prop:divnot0reverse} There exists a constant $C>1$, such
that, for any pair of cells $(K,L)$ and any function $\Psi:
{\Tr}^{N+1} \rightarrow {\mathbb R}_+$,
\begin{equation*}
\begin{split}
&C^{-1} \bigl[ 1- \|\div(a)\|_{\infty} \Delta t \bigr]^{N} {\mathbb
E}^{\mathbb Q}_{L} \bigl[ \Psi(\Gamma_0,\dots, \Gamma_N) {\mathbf
1}_{\{\Gamma_N=K\}} \bigr]
\\
&\hspace{5pt} \leq {\mathbb E}_{K} \bigl[ \Psi(\lKN_0,\dots, \lKN_N)
{\mathbf 1}_{\{K_N=L\}} \bigr] \\
&\hspace{10pt}
\leq C \bigl[ 1+ \|\div(a)\|_{\infty} \Delta t \bigr]^{N} {\mathbb
E}^{\mathbb Q}_{L} \bigl[ \Psi(\Gamma_0,\dots, \Gamma_N) {\mathbf
1}_{\{\Gamma_N=K\}} \bigr].
\end{split}
\end{equation*}
\end{prop}

Of course, Proposition \ref{prop:divnot0reverse} is weaker than
Proposition \ref{prop:invmeastorus}. Above, we are just able to
compare the law of $(\Gamma_0,\dots,\Gamma_N)$ under the measure
${\mathbf 1}_{\{\Gamma_N=K\}}\cdot{\mathbb Q}_L$ with the law of
$(K_N,\dots,K_0)$ under the measure ${\mathbf 1}_{\{K_N=L\}} \cdot
{\mathbb P}_K$. They are equivalent and the resulting density is
bounded from above and from below by positive deterministic
constants.
\vspace{5pt}
\\
\noindent {\bf Proof.} Without loss of generality, we can assume
that $\Psi={\mathbf 1}_{(J_0,\dots,J_N)}$, with $(J_0,\dots,J_N) \in
{\Tr}^{N+1}$, $J_0=L$ and $J_N=K$. Then,
\begin{equation}
\label{eq:reversedivnot0}
\begin{split}
{\mathbb E}_K \bigl[ \Psi(\lKN_0,\dots,\lKN_N) {\mathbf
1}_{\{K_N=L\}} \bigr] &={\mathbb P}_{J_N}
\{\lKN_{0}=J_0,\dots,\lKN_{N-1}=J_{N-1},\lKN_N=J_N\}
\\
&= {\mathbb P}_{J_N} \{K_0=J_N,K_1=J_{N-1},\dots,K_N=J_0\}
\\
&= p_{J_N,J_{N-1}} p_{J_{N-1},J_{N-2}} \dots p_{J_1,J_0}.
\end{split}
\end{equation}
By Proposition \ref{prop:-adivnot0} and \eqref{eq:reverse100}, we
know that
\begin{equation*}
\gamma_{J_n,J_{n+1}} = \frac{|J_{n+1}|}{|J_n|(1+\delta_{J_n} \Delta
t)} p_{J_{n+1},J_n}, \ 0 \leq n\leq N-1,
\end{equation*}
even if $J_n=J_{n+1}$.
 Plugging this relationship in
\eqref{eq:reversedivnot0}, we deduce that
\begin{equation*}
\begin{split}
{\mathbb E}_K \bigl[ \Psi(\lKN_0,\dots,\lKN_N) {\mathbf
1}_{\{K_N=L\}} \bigr] &= \bigl[ \prod_{n=0}^{N-1} (1 + \delta_{J_n}
\Delta t) \bigr] \frac{|J_0|}{|J_N|} \gamma_{J_0,J_1} \dots
\gamma_{J_{N-1},J_N}
\\
&= \bigl[ \prod_{n=0}^{N-1} (1 + \delta_{J_n} \Delta t) \bigr]
\frac{|J_0|}{|J_N|} {\mathbb E}^{\mathbb Q}_L \bigl[
\Psi(\Gamma_0,\dots,\Gamma_N) {\mathbf 1}_{\{\Gamma_N=K\}} \bigr].
\end{split}
\end{equation*}
By \eqref{highd:deltat} and \eqref{ass:iso} (recall that
\eqref{ass:iso} implies $\beta^{-1} h^d \leq |J| \leq h^d$ for all $J
\in {\Tr}$ and for some $\beta > 1$), we complete the proof.
\qed

\subsection{Analysis of the Gap}

Using the previous subsection, we analyze the term $R_N$ in
\eqref{eq:R+S}. For $N \geq 1$, we emphasize that
\begin{equation*}
\begin{split}
\sup_{1 \leq n\leq N} |R_n| &= \sup_{1 \leq n \leq N} \bigl|
\sum_{i=1}^{n-1} \bigl[ (e_{K_i} - x_{K_{i-1},K_i}) {\mathbf
1}_{\{K_{i-1} \not = K_i\}} \bigr] \bigr|
\\
&\leq 2 \sup_{1\leq n \leq N} \bigl|\sum_{i=n}^N \bigl[ ( e_{K_i} -
x_{K_{i-1},K_i}) {\mathbf 1}_{\{K_{i-1} \not = K_i\}} \bigr]\bigr|
\\
&= 2 \sup_{1 \leq n \leq N} \bigl|\sum_{i=0}^{n-1} \bigl[
(e_{\lKN_i} - x_{\lKN_{i},\lKN_{i+1}}) {\mathbf 1}_{\{\lKN_{i+1}
\not = \lKN_i\}} \bigr]\bigr|.
\end{split}
\end{equation*}
Keeping Proposition \ref{prop:divnot0reverse} in mind, we define the
analogue, but for the chain $(\Gamma_n)_{n \geq 0}$, that is
\begin{equation*}
\Xi_0=0, \ \Xi_n = \sum_{i=0}^{n-1} \bigl[ ( e_{\Gamma_i} -
x_{\Gamma_{i},\Gamma_{i+1}}) {\mathbf 1}_{\{\Gamma_{i+1} \not =
\Gamma_i\}} \bigr], \ n \geq 1,
\end{equation*}
By Proposition \ref{prop:divnot0reverse}, we deduce that, for any
starting cell $K$, for any terminal cell $L$ and for any $k \geq 0$,
\begin{equation}
\label{eq:gap} {\mathbb P}_K \bigl\{ \sup_{1 \leq n \leq N} |R_n|
\geq \frac{k}{2} h^{1/2}, K_N = L \bigr\} \leq C(1+\Delta t)^N
{\mathbb Q}_L \bigl\{
 \sup_{1 \leq n \leq N} \bigl|\Xi_n\bigr| \geq \frac{k}{4}
h^{1/2},\Gamma_N=K\bigr\}.
\end{equation}
Here is the main argument of the analysis.
\begin{prop} \label{prop:egamma} For a given cell $L \in
{\Tr}$, the process $(\Xi_n)_{n \geq 0}$ is a martingale
under ${\mathbb Q}_L$ with respect to the filtration
$(\sigma(\Gamma_0,\dots,\Gamma_n))_{n \geq 0}$.
\end{prop}
\noindent {\bf Proof.} The proof is quite obvious. For each $n \geq
0$, $\Xi_n$ is measurable with respect to
$\sigma(\Gamma_0,\dots,\Gamma_n)$ and
\begin{equation*}
\begin{split}
{\mathbb E}^{\mathbb Q}_L \bigl[ e_{\Gamma_n} -
x_{\Gamma_{n+1},\Gamma_n} | \Gamma_0,\dots,\Gamma_n \bigr] &=
\sum_{J \in \Gamma_n^+} \gamma_{\Gamma_n,J} ( e_{\Gamma_n} -
x_{J,\Gamma_n} )
\\
&=(1+\delta_{\Gamma_n} \Delta t)^{-1} \sum_{J \in \Gamma_n^+}
q_{\Gamma_n,J} ( e_{\Gamma_n} - x_{J,\Gamma_n} ) = 0. \qed
\end{split}
\end{equation*}

Following the proof of Lemma \ref{lem:S} and making use of
\eqref{highd:deltat} to bound the conditional variances of the
increments of $(\Xi_n)_{n \geq 0}$, we deduce
\begin{lem}
\label{prop:estimq} There exists a constant $C>0$ such that for any
$k \geq 0$
\begin{equation*}
\sup_{L} {\mathbb Q}_L \bigl\{ \sup_{0 \leq n \leq N} |\Xi_n| >
\frac{k}{4} h^{1/2} \bigr\} \leq C \bigl[ \exp \bigl( - \frac{k^2}{C N
\Delta t} \bigr) + \exp \bigl( - \frac{k}{C h^{1/2}} \bigr) \bigr].
\end{equation*}
\end{lem}

Gathering \eqref{eq:gap} and Lemmas \ref{lem:general40} and
\ref{prop:estimq}, we deduce (by modifying if necessary the constant
$C$ from line to line)
\begin{equation*}
\begin{split}
&\sum_{k \geq 0} \sum_{K \in {\Tr}_0} \sum_{L \in \Tr_k^N} (k+1)^p
{\mathbb P}_K \bigl\{ \sup_{0 \leq n \leq N} |R_n| \geq \frac{k}{2}
h^{1/2}, K_N =L\bigr\} \\
&\leq C (1+ \Delta t)^N \sum_{k \geq 0} \sum_{L \in \Tr_k^N}
\sum_{K \in {\Tr}_0} (k+1)^p {\mathbb Q}_L \bigl\{ \sup_{1 \leq n \leq
N} |\Xi_n| \geq \frac{k}{4} h^{1/2}, \Gamma_N =K\bigr\} \\
&\leq C (1+ \Delta t)^N \sum_{k \geq 0} \sharp[ {\Tr}^N_k]
(k+1)^p \sup_L {\mathbb Q}_L \bigl\{ \sup_{1 \leq n \leq N}
|\Xi_n| \geq \frac{k}{4} h^{1/2}\bigr\}
\\
&\leq C h^{-d/2} \exp(C N \Delta t) \sum_{k \geq 0} (k+1)^{d+p}
\bigl[ \exp \bigl( - \frac{k^2}{C N \Delta t} \bigr) + \exp \bigl( -
\frac{k}{C h^{1/2}} \bigr) \bigr].
\end{split}
\end{equation*}

By \eqref{eq:R+S} and \eqref{eq:Ssup}, we deduce (up to a new value
of $C$)
\begin{equation*}
\begin{split}
{\mathcal Q}_p^N &\leq Ch^{(p+1)/2} N \Delta t (1+ N \Delta t)^p
\\
&\hspace{15pt} + C h^{p/2} \exp(C N \Delta t) \sum_{k \geq 0}
(k+1)^{d+p} \bigl[ \exp \bigl( - \frac{k^2}{C N \Delta t} \bigr) +
\exp \bigl( - \frac{k}{C h^{1/2}} \bigr) \bigr].
\end{split}
\end{equation*}
Keeping the inequality $h \leq 1$ in mind and comparing the 
left-hand side below to an integral, there exists a positive constant
$C_p$, only depending on $p$ and on the same parameters as $C$, such
that
\begin{equation*}
 \sum_{k \geq 0}
(k+1)^{d+p} \bigl[ \exp \bigl( - \frac{k^2}{C N \Delta t} \bigr) +
\exp \bigl( - \frac{k}{C h^{1/2}} \bigr) \bigr] \leq C_p (1+N \Delta
t)^{(p+d+1)/2}.
\end{equation*}

Finally,
\begin{theorem}
\label{thm:general1} Assume \eqref{ass:iso}, \eqref{highd:deltat},
$h \leq 1$ and $\Delta t \leq \theta h$. Then, for any $p\geq 1$,
there exists a constant $C_p >0$, only depending on
$\|a\|_{\infty}$, $\alpha$, $d$, $\kappa$, $\eta$, $\theta$ and $p$,
such that, for all $N \geq 1$,
\begin{equation*}
h^{d/2} \sum_{K \in {\Tr}_0} {\mathbb E}_{K} \bigl[
\bigl|X_N - X_0 + \Delta t \sum_{n=0}^{N-1} a(X_n) \bigr|^p \bigr]
\leq C_p h^{p/2} \exp ( C_p N \Delta t ).
\end{equation*}
\end{theorem}
\noindent \vspace{5pt}

\subsection{Analysis of the Numerical Scheme}

We now prove the main results of the paper.

\begin{prop}
\label{prop:L_ploc} In addition to the assumptions of Theorem
\ref{thm:general1}, assume that $u^0$ is $\kappa$-Lipschitz
continuous. Then, for any $p\geq 1$, there exists a constant $C_p
>0$, only depending on $\|a\|_{\infty}$, $\alpha$, $d$, $\kappa$,
$\eta$, $\theta$ and $p$, such that, for all $N \geq 1$,
\begin{equation*}
\sup_{x \in \R^d} \bigl[ h^{-d/2} \sum_{K:|e_K-x| \leq h^{1/2}}
\|u^N_K - u(N \Delta t, \cdot)\|^p_{L^p(K)} \bigr]^{1/p} \leq C_p
h^{1/2} \exp(C_p N \Delta t).
\end{equation*}
\end{prop}
\noindent {\bf Proof.} By translation, it is sufficient to prove the
bound for $x=0$. To simplify the notations, we set for any random
variable $Y$ with values in $\R^d$:
\begin{equation*}
\|Y\|_{p,{\Tr}_0} = \bigl[ h^{d/2}
 \sum_{K \in {\Tr}_0} {\mathbb E}_{K} \bigl[ |Y|^p \bigr]
\bigr]^{1/p}.
\end{equation*}
(Above, 
${\mathcal T}_0$ is the set of cells $K$ such that $|e_K| \leq
h^{1/2}$.) 
Of course, $\| \cdot \|_{p,{\Tr}_0}$ is a norm on the space
of $\R^d$-valued random variables with a finite moment of order $p$
under every ${\mathbb P}_K$, $K \in {\Tr}_0$. By Theorem
\ref{thm:general1} and by the properties $h \leq 1$ and $\Delta t \leq
\theta h$, we have
\begin{equation*}
\forall N \geq 1, \ \bigl\| X_N - Z(N \Delta t,X_0) + \Delta t
\sum_{n=0}^{N-1} \bigl[a(X_n) - a(Z(n \Delta t,X_0)) \bigr]
\bigr\|_{p,{\Tr}_0}
 \leq C_p h^{1/2} \exp \bigl( C_p N \Delta t \bigr),
\end{equation*}
up to a new value of $C_p$. Following the proof of Proposition
\ref{prop:d1_general}, Gronwall's lemma yields
\begin{equation}
\label{eq:prout}
\forall N\geq 1, \ \bigl\| X_N - Z(N \Delta t,X_0)
\bigr\|_{p,{\Tr}_0} \leq C_p h^{1/2} \exp \bigl( C_p N \Delta
t \bigr),
\end{equation}
again for a new value of the constant $C_p$. Following the proof of
Theorem \ref{cor:L1conv},
\begin{equation*}
\begin{split}
\bigl[ h^{-d/2} \sum_{K \in {\Tr}_0} \|u^N_K - u(N
\Delta t, \cdot)\|^p_{L^p(K)} \bigr]^{1/p} \leq \kappa \bigl\| X_N -
Z(N \Delta t,X_0) \bigr\|_{p,{\Tr}_0} + C \exp(C N \Delta t) h,
\end{split}
\end{equation*}
for some $C>0$. This completes the proof. \qed \vspace{5pt}

As a by-product, we obtain the $L^{\infty}$ estimate announced in
Introduction: 
\begin{theorem}
\label{prop:L_inf}In addition to the assumptions of Theorem
\ref{thm:general1}, assume that $u^0$ is $\kappa$-Lipschitz
continuous. Then, for any $p\geq 1$, there exists a constant $C_p
>0$, only depending on $\|a\|_{\infty}$, $\alpha$, $d$, $\kappa$,
$\eta$, $\theta$ and $p$, such that, for all $N \geq 1$,
\begin{equation*}
\sup_{K \in {\Tr}} \sup_{x \in K} |u^N_K - u(N \Delta t, x)|
\leq C_p h^{(1-1/p)/2} \exp(C_p N \Delta t).
\end{equation*}
\end{theorem}
\noindent {\bf Proof.} Consider a given cell $K$. By Proposition
\ref{prop:L_ploc},
\begin{equation*}
h^{d/2p}\inf_{y \in K} |u^N_K - u(N \Delta t, y)| \leq
\bigl[ \beta h^{-d/2} \|u^N_K - u(N \Delta t, \cdot)\|^p_{L^p(K)}
\bigr]^{1/p} \leq C_p h^{1/2} \exp(C_p N \Delta t).
\end{equation*}
Since $u(N \Delta t,\cdot)$ is Lipschitz continuous with $\kappa
\exp(\kappa N \Delta t)$ as Lipschitz constant, we complete the
proof. \qed

By a regularization argument, we manage to weaken the required
assumption on $u^0$ in Proposition \ref{prop:L_ploc}:
\begin{theorem}
\label{th:L_1BV}
In addition to the assumptions of Theorem
\ref{thm:general1}, assume that $u^0$ belongs to $BV({\mathbb R}^d)$. 
Then, there exists a constant $C
>0$, only depending on $\|a\|_{\infty}$, $\alpha$, $d$, $\kappa$,
$\eta$, $\theta$ and the $BV$ semi-norm of $u^0$, such that, for all
$N \geq 1$, 
\begin{equation*}
\sum_{K \in {\mathcal T}}
\|u^N_K - u(N \Delta t, \cdot )\|_{L^1(K)} \leq C h^{1/2} \exp(C N
\Delta t). 
\end{equation*}
\end{theorem}

\noindent {\bf Proof.} 
The parameters $h$ and $N$ are fixed for the whole proof (with $h \leq
1$). The constants ``$C$'' and ``$C'$'' appearing below may depend on
the $BV$ semi-norm of $u^0$. By \cite[Chapter 5]{Ziemer}, we know that
the semi-norm $BV$, denoted by $\| \cdot \|_{BV(\R^d)}$, decreases by
convolution. In particular, we can find a smooth function $u_h^0 :
\R^d \rightarrow \R$, such that $\|u^0-u_h^0\|_{L^1(\R^d)} \leq
h^{1/2}$ and $\|\nabla u^0_h\|_{L^1(\R^d)} \leq \|u^0\|_{BV(\R^d)}$.
We then set, for all $x \in \R^d$, $\bar{u}^0(x) = |B(0,h^{1/2})|^{-1}
\int_{B(0,h^{1/2})} u_h^0(x-y) dy$. We let the reader check that
\begin{equation}
\label{final:1}
 \|u^0 - \bar{u}^0 \|_{L^1(\R^d)} \leq h^{1/2} \bigl( 1 +
\|u^0\|_{BV({\mathbb R}^d)} \bigr), \ \|\nabla \bar{u}^0\|_{L^1(\R^d)}
\leq \|u^0\|_{BV(\R^d)}.
\end{equation}
The reason why we introduce $\bar{u}^0$ is the following: the gradient
of $\bar{u}^0$ may be locally bounded, in the $L^{\infty}$ sense, by
local $L^1$ norms of $\nabla u_h^0$. Indeed, for any $x \in \R^d$ and
$R>0$, 
\begin{equation}
\label{final:1.1}
\sup_{y \in B(x,R)} |\nabla \bar{u}^0(y)| \leq |B(0,h^{1/2})|^{-1}
\|\nabla u_h^0\|_{L^1(B(x,R+h^{1/2}))}. 
\end{equation}

We denote by $\bar{u}(t,x) = \bar{u}^0(Z(t,x))$ the solution, at
$(t,x)$, of the transport problem with $\bar{u}^0$ as initial
condition and by $\bar{u}^N_K = {\mathbb E}_K(\bar{u}^0(K_N))$ the
corresponding approximate solution at time step $N$ in cell $K$. It is
clear that there exists a constant $C>0$ such that
\begin{equation}
\label{final:2}
\forall t \geq 0, \ \|\bar{u}(t,\cdot) - u(t,\cdot)\|_{L^1(\R^d)} \leq
\exp(C t) \|\bar{u}^0 - u^0 \|_{L^1(\R^d)}. 
\end{equation}
Similarly, with Proposition \ref{prop:-adivnot0} at hand,
\begin{equation*}
\begin{split}
\sum_{K \in {\mathcal T}} |K| |u^1_K- \bar{u}^1_K|
&\leq 
\sum_{K \in {\mathcal T}} |K| \sum_{L \in {\mathcal T}} p_{K,L} |u^0_L
- \bar{u}^0_L| \\
&= \sum_{L \in {\mathcal T}} |L| (1+\delta_L \Delta t) \sum_{K \in
{\mathcal T}} \gamma_{L,K} |u^0_L - \bar{u}^0_L| \leq (1+ C \Delta t)
\sum_{L \in {\mathcal T}} |L| |u^0_L - \bar{u}^0_L|. 
\end{split}
\end{equation*}
Iterating the procedure, we obtain
\begin{equation}
\label{final:3}
\sum_{K \in {\mathcal T}} |K| |u^N_K- \bar{u}^N_K|
\leq (1+ C \Delta t)^N \sum_{K \in {\mathcal T}} |K| |u^0_K -
\bar{u}^0_K| \leq (1+ C \Delta t)^N \|u^0 - \bar{u}^0 \|_{L^1(\R^d)}.
\end{equation}
By \eqref{final:1}, \eqref{final:2} and \eqref{final:3}, it is
sufficient to investigate $\sum_{K \in {\mathcal T}}
\|\bar{u}^N_K - \bar{u}(N \Delta t, \cdot) \|_{L^1(K)}$. For a cell $K$, 
we obtain by \eqref{final:1.1} and by stability of the solutions to
\eqref{carac}, for all $y \in K$,
\begin{equation*}
\begin{split}
|\bar{u}(N \Delta t,y) - \bar{u}(N \Delta t,e_K)|
&= |\bar{u}^0(Z(N \Delta t,y)) - \bar{u}^0(Z(N \Delta t,e_K))| \\ 
&\leq C \exp(C N \Delta t) h \sup_{|z-Z(N \Delta t,e_K)| \leq C
\exp(CN \Delta t) h} |\nabla \bar{u}^0(z)| \\
&\leq C \exp(C N \Delta t) h^{1-d/2} \|\nabla u_h^0 \|_{L^1(B(Z(N \Delta
t,e_K),C \exp(C N \Delta t) h^{1/2}))}, 
\end{split}
\end{equation*}
for some $C>0$ (which may vary from line to  line). Integrating with
respect to $y$, we have by inversion of the flow $Z(N \Delta t,
\cdot)$ 
\begin{equation*}
\begin{split}
&\sum_{K \in {\mathcal T}} \|\bar{u}(N \Delta t,\cdot) - \bar{u}(N
\Delta t,e_K)\|_{L^1(K)} \\
&\leq C \exp(C N \Delta t) h^{1+d/2} \int_{\R^d} \bigl[ |\nabla
u^0_h(z)| \sum_{K \in {\mathcal T}}  {\mathbf 1}_{\{|z-Z(N \Delta
t,e_K)| \leq C \exp(C N \Delta t) h^{1/2}\}} \bigr] dz \\
&\leq C' \exp(C' N \Delta t) h \|\nabla u_h^0 \|_{L^1(\R^d)}. 
\end{split}
\end{equation*}
It is thus sufficient to analyse  $\sum_{K \in {\mathcal T}} |K|
|\bar{u}^N_K - \bar{u}(N \Delta t,e_K)|$.
In what follows, we fix a point $x \in \R^d$ and we consider a cell $K
\in {\mathcal T}_x$, the set of cells $L$ such that $|e_L-x| \leq
h^{1/2}$. The triangular inequality yields (remind that $X_0=e_K$
under ${\mathbb P}_K$) 
\begin{equation}
\label{final:4}
\begin{split}
&|\bar{u}^N_K - \bar{u}(N \Delta t,e_K)| 
\\
&\hspace{15pt}
\leq \sum_{k \geq 0} {\mathbb
E}_K \bigl[ |\bar{u}^0_{K_N} - \bar{u}^0(Z(N \Delta t,X_0))| {\mathbf
1}_{\{k h^{1/2} \leq |X_N- Z(N \Delta t,X_0)| < (k+1) h^{1/2}\}}
\bigr]. 
\end{split}
\end{equation}
We are to bound, under ${\mathbb P}_K$, the difference
$|\bar{u}^0_{K_N} - \bar{u}^0(Z(N \Delta t,X_0))|$ on the set $\{k
h^{1/2} \leq |X_N- Z(N \Delta t,X_0)| < (k+1) h^{1/2}\}$ by the
gradient of $\bar{u}^0$ and by $(k+2) h^{1/2}$. (On this set, every
point in $K_N$ is at distance less than $(k+1)h^{1/2} +h \leq (k+2)
h^{1/2}$ from $Z(N \Delta t,X_0)$.) Since $X_0=e_K$ under ${\mathbb
P}_K$, we know by stability of the solutions to \eqref{carac}
that $|Z(N \Delta t,X_0) - Z(N \Delta t,x)| \leq C \exp(C N \Delta t)
h^{1/2}$ for some constant $C>0$. Therefore, up to a modification of
$C$, the ball $B(Z(N \Delta t,X_0),(k+3)h^{1/2})$ is included in the
ball $B_x^k = B(Z(N \Delta t,x),C\exp(CN \Delta t)(k+1)h^{1/2})$. By
\eqref{final:1.1}, we obtain (under ${\mathbb P}_K$)
\begin{equation}
\label{final:5}
\sup_{|z-Z(N \Delta t,X_0)| \leq (k+2)h^{1/2}} |\nabla \bar{u}^0(z)|
\leq C h^{-d/2} \|\nabla u_h^0\|_{L^1(B^k_x)}. 
\end{equation}
By \eqref{final:4} and \eqref{final:5},
\begin{equation*}
\begin{split}
&\sum_{K \in {\mathcal T}_x} |K| |\bar{u}^N_K - \bar{u}(N \Delta t,
e_K)| \\
&\hspace{15pt} \leq 
C h^{(1+d)/2} \sum_{k \geq 0} \biggl[ (k+2) \|\nabla
u_h^0\|_{L^1(B_x^k)} \sum_{K \in {\mathcal T}_x} {\mathbb P}_K \bigl\{
|X_N- Z(N \Delta t,X_0)| \geq k h^{1/2} \bigr\} \biggr]. 
\end{split}
\end{equation*}
Recall the Markov inequality: for a nonnegative random variable $Y$
and two reals $a,p>0$, ${\mathbb P}_K\{Y>a\} \leq a^{-p} {\mathbb
E}_K[Y^p]$. Choosing $Y=h^{-1/2} |X_N- Z(N \Delta t,X_0)| $, $a=k$ and
$p=d+3$, and referring to the proof of Proposition \ref{prop:L_ploc}
(see \eqref{eq:prout}), we obtain for $k \geq 1$: 
\begin{equation*}
\begin{split}
\sum_{K \in {\mathcal T}_x} {\mathbb P}_K \bigl\{ |X_N- Z(N \Delta
t,X_0)| \geq k h^{1/2} \bigr\} &\leq k^{-(d+3)} \sum_{K \in {\mathcal
T}_x} {\mathbb E}_K \bigl[ \bigl( h^{-1/2} |X_N- Z(N \Delta t,X_0)|
\bigr)^{d+3} \bigr] \\
&= k^{-(d+3)} h^{-d/2} \| h^{-1/2} (X_N- Z(N \Delta t,X_0))\,
\|_{d+3,{\mathcal T}_x}^{d+3} \\
&\leq C (k+1)^{-(d+3)} h^{-d/2} \exp(CN \Delta t).
\end{split}
\end{equation*}
The modification of $k$ into $k+1$ in the last line permits to recover
the case $k=0$. Finally,
\begin{equation*}
\sum_{K \in {\mathcal T}_x} |K| |\bar{u}^N_K - \bar{u}(N \Delta t,
e_K)| \leq C h^{1/2} \exp(C N \Delta t) \sum_{k \geq 0}
\bigl[ (k+2) (k+1)^{-(d+3)} \|\nabla u_h^0\|_{L^1(B_x^k)} \bigr]. 
\end{equation*}
When integrating the left-hand side with respect to $x$ over $\R^d$,
we obtain by Fubini's Theorem a term equal to a constant times
$h^{d/2} \sum_{K \in {\mathcal T}} |K| |\bar{u}^N_K - \bar{u}(N \Delta t,
e_K)|$. When integrating the right-hand side, we have by
inversion of the flow $Z(N \Delta t,\cdot)$ 
\begin{equation*}
\begin{split}
\int_{\R^d} \|\nabla u_h^0\|_{L^1(B_x^k)} dx &= \int_{\R^d} \int_{B(Z
(N \Delta t,x),C\exp(CN\Delta t)(k+1)h^{1/2})} |\nabla u_h^0|(z) dz dx
\\ 
&\leq C' \exp(C' N \Delta t) \int_{\R^d}
\int_{B(x,C\exp(CN\Delta t)(k+1)h^{1/2})}|\nabla u_h^0|(z) dz dx \\
&\leq C' (k+1)^d h^{d/2} \exp(C' N \Delta t) \|\nabla
u_h^0\|_{L^1(\R^d)}. 
\end{split}
\end{equation*}
This completes the proof. \qed

\section{Annex}\label{annex}

We now show the concentration inequality used to prove Lemmas
\ref{lem:S} and \ref{prop:estimq}.

\begin{prop}
\label{prop:concentration} Let $(\Omega,{\mathcal A},{\mathbb P})$
be a probability space and $d \in \N \setminus \{0\}$. Then,
there exists a constant $c>0$, only depending on $d$, such that, for
any ${\mathbb R}^d$-valued martingale $(Y_n)_{n \geq 0}$ with
respect to a given filtration $({\mathcal H}_n)_{n \geq 0}$,
satisfying $Y_0=0$ and, for all $n \geq 0$, $|Y_{n+1}-Y_n| \leq 1$
and ${\mathbb E}[ |Y_{n+1}-Y_n|^2|{\mathcal H}_n] \leq v$ for a
deterministic real $v
>0$, the following holds for all $n \geq 0$:
\begin{equation*}
\forall u >0, \ {\mathbb P} \bigl\{ \sup_{0 \leq k \leq n} |Y_k|
\geq u \bigr\} \leq c \bigl[ \exp\bigl( - \frac{u^2}{c n v} \bigl) +
\exp \bigl( - \frac{u}{c} \bigr) \bigr].
\end{equation*}
\end{prop}
\noindent{\bf Proof.} Without loss of generality, we can assume that
$d=1$. Indeed,
\begin{equation*}
\forall u >0, \ {\mathbb P} \bigl\{ \sup_{0 \leq k \leq n} |Y_k|
\geq u \bigr\} \leq \sum_{i=1}^d {\mathbb P} \bigl\{ \sup_{0 \leq k
\leq n} |(Y_k)_i| \geq u d^{-1/2} \bigr\}.
\end{equation*}
In the one-dimensional case, it is sufficient to investigate
${\mathbb P} \{ \sup_{0 \leq k \leq n} Y_k \geq u \}$, the lower bound
following from an obvious change of sign. We then apply
\cite[Proposition 1.6]{freedman} (choose $X_k=0$ for $k>n$, $a=u$ and $b=nv$
 in the statement 
of \cite[Proposition 1.6]{freedman})
\begin{equation*}
{\mathbb P} \bigl\{ \sup_{0 \leq k \leq n} Y_k \geq u \bigr\} \leq
\exp \bigl( - \frac{u^2}{2(u + n v)} \bigr). 
\end{equation*}
There are two cases: either $u \leq nv$ or $u > nv$. We obtain
\begin{equation*}
{\mathbb P} \bigl\{ \sup_{0 \leq k \leq n} Y_k \geq u \bigr\} \leq
\exp \bigl( - \frac{u^2}{4 n v} \bigr) + \exp \bigl( - \frac{u}{4}
\bigr). 
\end{equation*}
This completes the proof. \qed


\end{document}